# The Anthyphairetic Revolutions of the Platonic Ideas

**Stelios Negrepontis**

θεὰ σκέδασ' ἠέρα, εἴσατο δὲ χθών·
*Odusseia*, Book XIII, line 352

**Abstract.** In the present work Plato's philosophy is interpreted as an imitation, a close philosophic analogue of the geometric concept of a pair of lines incommensurable in length only and of its (palindromically) periodic anthyphairesis. It is shown, by an examination of the Platonic dialogues *Theaetetus*, *Sophistes*, *Politicus*, and *Philebus*, that
(a) a Platonic Idea is the philosophic analogue of a dyad of lines incommensurable in length only,
(b) the Division and Collection, the method by which humans obtain knowledge of a Platonic Idea, is the philosophic analogue of the palindromically periodic anthyphairesis of this dyad, and
(c) a Platonic Idea is One in the sense of the self-similarity induced by periodic anthyphairesis.

A byproduct of the above analysis is that
(d) Theaetetus had obtained a proof of the Proposition: The anthyphairesis of a dyad of lines incommensurable in length only is palindromically periodic.
It is further verified that the concepts and tools contained in the Theaetetean Book X of the *Elements* suffice for the proof of the Proposition.

**Outline**. According to the *Philebus* 16c a **Platonic Idea** is the **mixture** of the two principles **Infinite** ('apeiron') and the **Finite** ('peras') (Section 2). A key step to this interpretation is the discovery, that, according to the *Philebus* 23b-25e, these two principles are close philosophic analogues of the concepts of **finite and infinite anthyphairesis (commensurability and incommensurability**, accordingly), described in Propositions 1-8 of Book X in Euclid's *Elements* (Section 3).

Since the components that enter in the mixing have so heavily a mathematical content, it is natural to try to locate first a mathematical mixture of commensurability and incommensurability. Set in this way, the answer lies in front of us in the concept of **incommensurability in length only of two lines** a,b (Definitions 2, 3 in Book X of the *Elements*), a mixture of incommensurability (a, b) and commensurability ($a^2$, $b^2$) (Section 4).

The main bulk of the argument (Sections 5-12), consists in showing that a Platonic Idea is a mixture of the Infinite and the Finite in the sense that it is the philosophic analogue of a dyad of lines incommensurable in length only.

In the *Theaetetus* 147d-148d, where the mixture, called **'dunamis'**, is introduced by Theaetetus, Plato gives an account of the incommensurabity proofs, initially by Theodorus and then by Theaetetus, of a, b in case $a^2=Nb^2$ for non-square N; but the mathematical result achieved by Theaetetus is couched in unclear philosophical terms, involving the crucial phrase on the collection into One of powers (‘hai dunameis’, 147d7-8) infinite in multitude. Delicate arguments, that have totally escaped the attention of the scholars, show that the **plural ‘hai dunameis’**, employed in this description, is **distributive,** and that Theaetetus succeeded in **Collecting into One every power separately** that appeared **infinite in multitude** by Theodorus’ proof of incommensurability (Section 5).

Socrates, at the end of the *Theaetetus* passage, asks that Theaetetus **imitate (‘peiro mimoumenos’**, 148d4) his method in order to tackle the philosophic problem of knowledge of the Ideas (‘episteme’), and the imitation, the method of Division and Collection, by which ‘episteme’ is achieved, is presented in the sequel, *Sophistes* and *Politicus*, of the Platonic trilogy (Section 6).

The distributive interpretation of the *Theaetetus* 147d7-8, the imitation of the resulting mathematical Division and Collection to produce the philosophical one, and the Platonic accounts of the philosophical Division and Collection lead to three **plausibility arguments** suggesting that Theaetetus mathematical achievement consisted in proving that every dyad of incommensurable in line only lines possess a **periodic anthyphairesis**, These plausibility arguments take us within striking distance from the modern theorem (established by the successive efforts of great mathematicians (including Fermat, Brouncker, Euler, Lagrange, Legendre, Galois): every pair of lines incommensurable in length only possesses **(palindromically)** periodic anthyphairesis (Section 7).

How would Theaetetus prove periodicity in anthyphairesis? Reasonable arguments lead to the conclusion that, unlike the Pythagoreans and unlike Theodorus, Theaetetus had created a theory of ratios of magnitudes (exactly the one reported in the *Topics* 158b-159a), based on equal anthyphairesis. In such a theory the most natural way of proving periodicity is by means of the **Logos Criterion** (Section 8).

The most decisive arguments, confirming the plausibility arguments developed in Section 7, are contained in the analytic examples of Division and Collection in the *Sophistes* and the *Politicus*. The knowledge of the **Angler** (Section 9) and of the **Sophist** (Section 10) in the *Sophistes*, paradigms of Division and Collection, are shown to be close philosophic analogues of anthyphairesis, achieving periodicity by the Logos criterion (the second one by employing the fundamental analogy of the divided line of the *Politeia* 509d-510b). But, and this is nothing short of astounding, in the *Politicus* the knowledge of the **Statesman** is achieved by a further refinement of the method, imitating not just periodicity but in fact **palindromic periodicity** (Section 11). Thus Plato makes sure that Theaetetus has proved the full palindromicity theorem, stated in Section 7. The **Double Measurement** in the *Politicus* 283a-287b philosophical mixture (even manages to take a form similar to the mathematical one (Section 12).

On the basis of the plausibility arguments, but mostly on the basis of the evidence provided by the divisions and Collections in the *Sophistes* and *Politicus*, we then have little choice but to conclude that Theaetetus had in fact proved the full palindromic periodicity theorem. The introduction of the **‘apotome’** (Proposition 73) and **‘binomial’** (Proposition 36) lines and their **conjugacy** (Propositions 112-114) in Book X. provide the necessary tools for a natural reconstruction (Section 13).

Several Platonists, in perplexed puzzlement, rejected the notion that a Platonic Idea could be a mixture of the infinite and the finite and divisible, in fact ad infinitum, since such properties appear to preclude the sine qua non of a Platonic Idea, its Oneness. What they have failed to understand is that an entity divisible, but in endless revolution and periodicity while being divided, is an entity possessing the most satisfying and internally determined concept of **Oneness, self-similarity** (Section 14).

## 1. The basics of Plato's philosophy

According to Plato's philosophy
(a) every intelligible Platonic idea is a One, a monad, always the same, not subject to generation or destruction (e.g. *Symposium* 210e-211b, *Phaedo* 78d,80b, *Philebus* 15b1-4);
(b) Ideas are knowable ('episteme') to humans with the method of Division and Collection (*Theaetetus*, *Sophistes*, *Politicus*, *Philebus*, *Phaedrus*, *Parmenides*); and
(c) The sensibles are subject to generation, destruction, and continued change. Sensibles participate in the intelligibles (e.g. *Phaedo* 100c-e).

It is the purpose of this paper is to examine the intelligible, Platonic Ideas, and to show that Ideas and their knowability must be understood as close philosophic analogues of commensurable in power only dyads of lines and their periodic anthyphairesis[1]. We are not concerned in the present work with a corresponding interpretation of the sensibles.[2]

---

[1] Relevant definitions and concepts are given in Section 6.
[2] The Platonic description of the sensibles (in the *Timaeus),* and their participation in the Platonic Ideas are to be found in the anthyphairetic approximations (the modern convegents of the continued fractions) of commensurability in power only (cf. Negrepontis [N, 2005] (Section F) for an outline).

## 2. Every Platonic Idea is the Philebean mixture of the infinite and the finite (*Philebus* 16c)

According to the *Philebus*, each of the Platonic Ideas is a **mixture** of the Philebean **infinite** and the Philebean **finite**.[3]

***Philebus* 16c**
***Socrates.*** *One which is easy to point out, but very difficult to follow for through it* ***all the inventions of art*** *have been brought to light. See this is the road I mean.*
***Protarchus.*** *Go on what is it?*
***Socrates.*** *A gift of gods to men, as I believe, was tossed down from some divine source through the agency of a Prometheus together with a gleaming fire; and the ancients, who were better than we and lived nearer the gods, handed down the tradition that* ***all the things which are ever said to exist*** *('ton aei legomenon einai') are* ***consist of one and many*** *and* ***have inherent in them the finite and the infinite ('peras de kai apeirian en autois sumphuton echonton').***

***Philebus* 23c-d**
*Socrates. We said that God revealed in the universe two elements,* ***the infinite and the finite****, did we not?*
*Protarchus. Certainly.*
*Socrates. Let us, then, assume these as two of our classes, and a third,* ***made One by combining ('hen ti summisgomenon') these two.***

***Philebus* 25b**
***Socrates.****Well, what shall we say is the nature of* ***the third class, made by combining ('to meikton') these two****?*

Thus, we will be able to discover the nature of Platonic Ideas and Beings, if we understand the meaning of the terms 'Philebean finite', 'Philebean infinite' and 'Philebean mixture' of infinite and finite.

---

[3] Some scholars, e.g. Cherniss [Ch, 1945], Ryle [Ry, 1966], do not think that 'ta aei legomena einai' refer to the Platonic Ideas; this position is to be understood within a general view that, Plato's repeated dithyrambic descriptions of Division and Collection not withstanding, downgrades the method. Cf. Sections 12. 5 and 14.4-5 below.

## 3. The Philebean infinite and finite is the philosophic analogue of the incommensurability and commensurability, respectively (*Philebus* 23b-25e)

During the fall semester of 1996, while being Visiting Professor at the Mathematics Department of Cyprus University, I discovered, in the Platonic dialogue *Philebus* 23b-25e, the first hard evidence that Plato's philosophy was related to the geometric concept of anthyphairesis, found in the *Elements* (finite Euclidean algorithm for numbers in Book VII, finite or infinite anthyphairesis of magnitudes in Book X). I gave a lecture on this subject at the Philosophy Department of Cyprus University on December 2**,** 1996**,** with the characteristic title:

*'A Mathematician Reads the* Philebus*'.*

An extended version of my lecture has appeared in [N, 1999].In the period 1996-2000 that followed, I wrote several versions of this discovery and gave a large number of lectures on it. My interpretation of the Philebean infinite and finite is adequately described in the two extracts I supply below (translated from the Greek) from two publications.

In the paper [N, 2000] (pp. 16-20), I write on the *Philebus* 23b-25e):

***'1a. The Dichotomy of Infinite and Finite in Euclid.***
*In the definitions 1 and 2 of Book X of the* Elements *a fundamental dichotomy s described. We can state this dichotomy as follows:*
*We consider as Whole the set R of all ordered pairs (a,b) of line segments… with a>b. The first fundamental dichotomy of R (given by definition 1) is between the class A of all incommensurable pairs and the class P of all commensurable pairs..*
*Propositions 1-8 of Book X are devoted to the proof of the following basic characterizations:*
*(a) P=the class…of all (a,b) in R with finite anthyphairesis=the class of all pairs (a,b) in R that are 'as number to number', and*
*(b) A= the class…of all (a,b) in R with infinite anthyphairesis=the class of all pairs (a,b) in R that are 'not as number to number'.*

***1b. The Dichotomy of Infinite and Finite in Plato****.*
*In passage 23b-25e of the Platonic dialogue* Philebus *there is a treatment of the Platonic dichotomy of infinite and finite.*
*The Whole, to which the the platonic dichotomy applies, consists of all ordered pairs of the type (more, less). Such a pair is infinite in case the relation (more, less) is renewed at every stage of the process and finite if the relation (more,less) stops or ends at some stage.*
*And it is finite, Plato informs us further, if each element of the pair (more, less) is a quantity (with respect to some measure) (cf. anonymous scholion 39 in Book X of the* Elements*), alternately, if the elements of the pair are as number to number, alternatively, if the pair is commensurable.*
*It becomes evident (after a somewhat more extended and detailed analysis into which we have no time to enter) that the Platonic Infinite (resp., the Platonic Finite) coincides with the Euclidean Infinite A (resp., the Euclidean Finite P)'*[4]

---

[4] In the newspaper *Avgi* there was an extensive presentation-interview of my interpretation of Plato's philosophy in three Sunday issues (July 30, August 6 and 13, 200) [Av, 2000]. In this presentation (July 30, 2000) I write on the *Philebus* 23b-25e:

In [N, 2005] I published a detailed account of my interpretation, where I have added two new items, in relation to the previous publications:
(a) I have included Proclus' *Commentary to Euclid* 5,11-7,12, confirming my interpretation of the Philebean principles of Finite and Infinite in terms of incommensurability and commensurability, respectively; and
(b) I have attempted the rather ambitious argument of showing directly, not through its opposition to the Finite, that Philebean Infinite is the philosophic analogue of infinite anthyphairesis.

For (a) I write in [N, 2005]:

*'Proclus' anthyphairetic account of the Philebean principles of the Finite and the Infinite…*
Full confirmation of our interpretation is supplied by Proclus, in his work *Comments to Euclid* **5,11-7,12,** where it is stated
(i) that the Philebean principle of the Finite is the cause of commensurability of magnitudes **(7,1-5),** and
(ii) that the Philebean principle of the Infinite is the cause of incommensurability of magnitudes **(6,19-21).**
These comments connect in a causal relation the Philebean philosophical Infinite with the mathematical concept of incommensurability, which we know to be equivalent, according to propositions X.2-3 of the *Elements*, with the mathematical concept of infinite anthyphairesis, so that the cause of mathematical infinite anthyphairesis is the Philebean principle of the Infinite; and also it connects in a similar relation the Philebean philosophical Finite with the mathematical concept of commensurability, which we know to be equivalent, according to propositions X.1-8 of the *Elements*, with the mathematical concept of finite anthyphairesis, so that the cause of mathematical finite anthyphairesis is the Philebean principle of the Finite.'

**Note on existing interpretations**. To the best of my knowledge the anthyphairetic interpretation of the *Philebus* 23b-25e is completely original, as no existing interpretation relates it with anthyphairesis. In fact, some of the classical interpretations, by Hackforth [Ha, 1945], G. Striker [St, 1970], Gosling [Go, 1975], D. Frede [FrD,1992], S. Delcomminette [De, 2006] (to name just a few), make no correlation to Euclid's *Elements*, Book X or to anthyphairesis; and, the recent volume, edited by J. Dillon and L. Brisson [DB, 2010], with contributions of some 52 papers by 42 scholars on various aspects of the *Philebus*, contains not even one direct or

---

'In Euclid we have a detailed description of the dichotomy between finite and infinite anthyphairesis… The infinite and the finite occur everywhere in Plato, but the most revealing text is undoubtedly the passage in the *Philebus* 23b-25e, where Plato reproduces in a very faithful manner, essentially the same terms that Euclid uses in Book X of the *Elements* in order to describe the fundamental dichotomy between infinite and finite anthyphairesis. The only object that may be described by the terms: finite, rest, end of every opposition, commensurable, each a quantity, as number to number, as measure to measure is the finite anthyphairesis'.

indirect reference to anthyphairesis or to Book X of the *Elements* or to Proclus' Comments on the Philebean infinite and finite.[5]

**3a**. In section 12.1 below we will make use of the following statement from the *Philebus*:

'Every pleasure ('hedone'), or more precisely every dyad (pleasure and pain) ('hedone kai lupe'), is an infinite'.

***Philebus* 27e**
Socrates. Is **pleasure and pain** a finite, or are they among the things which **admit of more and less**?
Philebus. Yes, they are among those which **admit of the more**, Socrates; for **pleasure** would not be absolute good if it were not **infinite in multitude and in the more**.

***Philebus* 31a**
Socrates**.** Let us, then, remember … that **pleasure** was itself **infinite** and belonged to the class which, in and by itself, has not and never will have either beginning or middle or end.

***Philebus* 41d**
Socrates. And have we not also said and agreed and settled something further?
Protarchus. What?
Socrates. That both **pleasure and pain admit of the more and less** and are of the class of **the infinite**.
Protarchus. Yes, we have said that, certainly.

[5] I presented my anthyphairetic interpretation of the *Philebus* 23b-25e (and of other Platonic passages), in the seminar on Philosophy at the National Technical University (Athens), co-directed by V. Karasmanis, on May 1999.
It was recently brought to my attention that V. Karasmanis has published a paper (entitled *Continuity and Incommensurability in Ancient Greek Philosophy and Mathematics*, in the volume edited by G. Anagnostopoulos, *Socratic, Platonic and Aristotelian Studies: Essays in Honor of Gerasimos Santas*, Philosophical Studies Series, volume 117, Springer 2011, p. 389-399), where he is (a) presenting an interpretation of the *Philebus* 23b-25e in terms of anthyphairesis, using the arguments I have already employed, but (b) making no reference to my original (1996) discovery, to my publications [N, 1999], [N, 2000], [N, 2005] or to published accounts [Av, 2000], and to my May 1999 lecture in his seminar on this interpretation.

## 4. The Theaetetean mixture of commensurable and incommensurable (*Theaetetus* 147e2-b3, *Elements*, def. 2, 3 in Book X)

After this identification, the question is what might be the Philebean mixture of the (Philebean) infinite and the (Philebean) finite. If we succeed in identifying the nature of the Philebean mixture, we will be able to understand the meaning of the Platonic Idea.

The two Philebean principles infinite and finite are the philosophic analogues, the imitations of the geometric concepts incommensurability and commensurability, correspondingly. It is then natural to conjecture that the philosophical Philebean mixture of the Philebean infinite and Philebean finite will be the philosophic analogue, the imitation of a suitable geometric mixture of commensurability and incommensurability.

Thus it is natural to look for a mixture of commensurability and incommensurability.

In definitions 2 and 3 of Book X of the *Elements* the central concept of Book X of the *Elements* is introduced, the **commensurability in power only** of two lines, equivalently the **incommensurability in length only** of two lines.

'(Two) straight-lines are commensurable in power ('dunamei summetroi') when the squares on them are measured by the same area, but (are) incommensurable (in power) when no area admits to be a common measure of the squares on them'.

**Incommensurable in length only**, equivalently **commensurable in power only**, are two line segments a and b such that a, b are incommensurable and their squares $a^2$, $b^2$ are commensurable.

Incommensurability in length only is indeed a geometric mixture of incommensurability and commensurability.

The concept is due to Theaetetus, as is not the case with the concepts commensurable/incommensurable, which are Pythagorean in origin. The circumstances of the introduction this concept by Theaetetus are related to us by Plato in the *Theaetetus* 147d3-148b4, a difficult and ill understood passage that nevertheless needs be clarified.

Indeed, Theaetetus introduced the concept of **power ('dunamis')**, by definition a line a, such that a is incommensurable in length ('mekei ou summetrous') to the one foot line b, but whose square $a^2$ is commensurable ('tois d' epipedois ha dunantai [summetrous]') to the square foot $b^2$ (147e2-148b4).

We will refer to a pair of lines a and b, such that a, b is incommensurable and $a^2$, $b^2$ commensurable as the **Theaetetean mixture** of commensurability and incommensurability.

## 5. The distributive plural argument: <u>Every</u> dunamis (Theaetetean mixture) possesses division and collection, <u>not</u> all powers collectively (*Theaetetus* 147d3-148b4)

In the *Theaetetus* 147d3-148b4 Plato relates to us not only the fact that Theaetetus introduced the notion of dunamis (Theaetetean mixture), but also the reason and the circumstances for doing so. During a lesson on some initial quadratic incommensurabilities taught by Theodorus, Theaetetus (and his companion) had a genial idea that that led him first to the introduction of the concept of power (dunamis, the Theaetetean mixture), and second to the proof of a mathematical result, relying on the concept of power, whose statement however, as given by Plato, involves the 'collection of the infinite multitude of powers into One', unfortunately a statement mathematically unclear.

Although the, admittedly difficult and delicate, passage in question is quite crucial, nevertheless it has been seriously misunderstood and compromised by most platonic scholars.
The first step in clarifying Theaetetus mathematical result is to decide on a question of linguistics: by the statement

'ai dunameis ephainonto apeiroi to plethos' (147d7-8)

does Plato want
(i) to say that the 'dunameis' apeeared to be apeiroi to plethos' as a totality, as a set, as it would appear more natural, or
(ii) to say that every power appeared to be infinite in multitude, as it would appear less natural (at least to our modern ears).

Technically we might say that the question is if the plural 'hai dunameis' is **distributive** (case (ii)),
or **collective** (case (i)).[6]
The passage has been interpreted collectively essentially by all scholars, including Knorr, Burnyeat. (Cf. Note in this Section below)
We will argue now that, on the contrary, Plato means to read the plural in a distributive manner; in which case, the *Theaetetus* 147d3-148b3 passage is to be rendered as follows:

**(147d3-6)** Theodorus in his lesson had shown ('apophainon') that if a ia line such that $a^2=Nb^2$, where N is one of the non-square numbers from N=3 to N=17 and b is the line of length one foot, then a is **incommensurable in length** ('mekei ou summetroi') to the one foot line;

**(147d7-e3)** Because, as a result of Theodorus lesson, it appeared that every power was divided into **an infinite multitude** of parts ('**apeiroi to plethos'**) (147d7-8), Theaetetus and his companion had the idea and in fact succeeded in **collecting** every power into One ('sullabein eis hen') (147d8-e3);

---

6 For example: 'the policemen wear uniform' is a collective plural; 'the policemen have surrounded the block' is a distributive plural.

**(147e3-148b3)** Theaetetus succeeded in the (Division and) Collection into One of every power by realizing that every **power ('dunamis')** is precicely a line segment incommensurable in length only to the one foot line, **the Theaetetean mixture** of incommensurability and commensurability.

Support for this distributive interpretation is based on the following three pieces of evidence, which have been ignored by all scholars so far [7]:

**[A]** The plural about powers is distributive
in the first occurrence ('[dunameis] mekei ou summetroi te podiaia', 147d4-5), and
in the third occurrence ('[dunameis]…mekei men ou summetrous ekeinais, tois d' epipedois ha dunantai [summetrous]', 148b1-2),
since every power is incommensurable or commensurable to the one foot, and not all the powers together;
hence it is natural to conclude that the plural about powers
in the second in between occurrence (hai dunameis apeiroi to plethos, 147d7-8)
is distributive as well.

**[B]** ***Scholia in Platonem [Theaetetus 162e6-7]*** **(SIP)**

Plato comments briefly on Theodorus' demonstrations in the following *Theaetetus* 162e6-7 passage:

*'Imagine how utterly worthless Theodorus or any geometer would be if he were prepared to rely on probability to do geometry.'*

With this passage (162e6-7), Plato clearly refers to the Theodorus' incommensurability demonstrations mentioned earlier in the *Theaetetus* 147d3-e1, the only Theodorus demonstrations mentioned in the *Theaetetus*. An interesting commentary of the 162e6-7, found in the anonymous work *Scholia in Platonem* (abbreviated SIP), runs as follows:

*'If we accept the judgement of the many as the dominant one in geometry, we would be ridiculous to claim that*
***[SIP1]*** *maginitudes are* ***incommensurable*** *to each other,*
*and that*
***[SIP2] the finite line is divisible ad infinitum*** *('diaireten eis apeiron'),*
*and* ***the like*** *('ta toiauta').'*

The following reading of the SIP commentary is suggested:
Theodorus, when proving that each of the powers (of 3,5,…,17 feet) and the one foot line (the magnitudes of SIP) are **incommensurable** to each other ([SIP1], as in statement 147d1-3)), employed,
not any considerations of probability (which would be worthless and make him ridiculous),
but logically impeccable proofs ('apodeixin de kai anagken', 162e4-5) involving the **division *ad infinitum*** of each of these finite lines ([SIP2],
so that every one of these lines, each power, is **infinite in multitude**, as in statement 147d7-8)).
The reading leads to a distributive interpretation of 'hai dunameis' in 147d7-8.

[7] In Section 6.3 below we will give a fourth argument in favor of distributivity of the plural.

**[C]** The *Anonymi Commentarius in Platonis Theaetetum* (ACIPT), possibly Eudorus, 1st century BC[8], lends strong support to the distributive interpretation. Indeed, commenting on the infinity arising from the 'powers infinite in multitude', and their collection achieved by Theaetetus and his companion, in the *Theaetetus* 147d-148b passage, he makes the following statements:
[1] an infinity in general must be collected, if possible;[9]
[2] a line admits of infinite either by (infinite) division, or by (infinite) increase;[10]
[3] the collection[11] of an infinity in a line can be achieved by passing over to numbers;[12]
[4] a line infinite by increase cannot be collected by numbers;[13]
[5] the collection of an infinity in a line can be achieved by passing over to numbers, because numbers are commensurable to each other;[14]
[6] the collection of an infinity in a line is achieved by passing from the less clear magnitudes to the more clear numbers[15], namely from incommensurables to commensurables, that is to 'as number to number'[16].

By statement [2] we conclude that
[7] the plural in the expression 'the powers infinite in multitude, as already suspected, is distributive and 'every power infinite in multitude' is meant;

by statement [4] and the fact that Theaetetus and companion succeeded in collecting the infinite in multitude powers, we conclude that
[8] every power is infinite in power by division *ad infinitum*;

by statements [3], [5], and [6] we conclude that
[9] the infinity of every power is due to incommensurability; and,

by statements [3], [5], [6], we conclude that
[10] 'collection into One' concerns the infinite multitude of parts of a power brought about by Theodorus incommensurability proofs, and is made possible because a power is not just any line incommensurable in length, but a line whose square is commensurable, namely a Theaetetean mixture of commensurability and incommensurability..

Even with the (novel) distributive reading of 'hai dunameis' in 147d7-8, still the meaning of the central term, 'collection into One', is not clear, as yet. But it is clear

---

8 cf. Tarrant [Tar, 1985]
9 ACIPT 37, 3-12
10 ACIPT 36,45-48 ('the lines admit of the indefinite either by increasing or by dividing them').
11 collection renders 'perilambanein' (as in the *Theaetetus* 148d6) (26,11; 37,5; 37,10; 37,29; 37,43; 37,46; 45,48; 46,39), periorizein (42,32) and horizein (37,1; 37,11).
12 ACIPT 36,48-37,3 ([and because the lines] are limited by numbers, Theaetetus and his companion passed over to numbers'); 42, 30-33
13 ACIPT 37, 39-44
14 ACIPT 26,13-18 ('thus they came to number as a consequence of the fact that all the numbers are commensurable to each other')
15 ACIPT 32,1-4; 44,50-45,3 ('from the less clear one must pass to the more clear as when passing from magnitudes to numbers'). Exactly the same point is made in the *Scholion* X.39 to the *Elements*.
16 ACIPT 41,40-42,18 (in analogy to the description for cubes: 'they came to the cubes', as 'the solid itself is commensurable to the solid; for it has ratio, as number to number, while th sides are incommensurable') that it has to do with a proof of the incommensurability of powers (after all the original lesson by Theodorus was about incommensurability of powers, and Proposition X.9 of the *Elements* is credited to Theaetetus). The fact that Theodorus' proof of incommensurability led him to a division ad infinitum certainly suggests an **anthyphairetic** approach to incommensurability, namely a proof employing Proposition X.2 of the *Elements*.[17] But we will not consider this evidence as yet conclusive, and we will not rely on it . The exact content of Theaetetus mathematical contribution will become clear later, in Sections 12 and 13.

**Note on existing interpretations**. To the best of my knowledge the interpretation of the *Theaetetus* 147d-148b passage, based on a distributive reading of the sentence 147d7-8 and based on arguments [A], [B], [C], above, is completely original. A preliminary version of this argument appeared in the Negrepontis-Tassopoulos paper [NT, 1212]. The crucial sentence is universally interpreted collectively, and not distributively, by all scholars.Thus

Cornford [Co,1935], p.23, translates: 'seeing that these square roots were evidently infinite in number'.

Klein [Kl, 1977}writes (p.80): 'Since the multitude of such incommensurable lines had appeared to Theaetetus and to young Socrates to be infinite, they had tried to find a way to encompass all of them in some unity,so that they could be called by one name.'

Sedley [Se, 2002] writes: 'Theaetetus recounts how he and his fellow student the younger Socrates were set by Theodorus the task of generalizing over the infinitely many square areas with integer areas but irrational sides.'

Tschemplik [Tsc, 2008], p.161, writes: 'Theaetetus already knew that he could not solve Socrates' questions the same way that he solved the mathematical questions posed to him byTheodorus (148e). This fact indicates an early awareness about the difference between mathematics and philosophy.Rather than looking at Theodorus's demonstrations and arriving at a solution to the problem of irrational numbers by formulating a hypothesis which comprehends all irrational numbers, dividing them into square and oblong numbers, Socrates is asking Theaetetus to examine his own activity, not just practiceit, when he is asking him to answer the question about knowledge.' (p.161)
It is interesting that the insistent exhortation for imitation of mathematics to philosophy is interpreted and perceived as a sign of difference between mathematics and philosophy.

Stern [Ste, 2008] writes (p. 62): 'Together, Theaetetus and young Socrates produce what Theodorus, who did not participate in the postlecture conversation, did not even attempt – specifically, a way of thinking of the infinite series of such numbers that gathers them into a unity. Theaetetus thus shows himself capable of performing the crucial act of understanding, gathering a many into a one.’
The statement makes very little sense. For a mathematician like Theaetetus to give a one name, dunamis, to all the incommensurable lines considered would be a triviality, not a significant mathematical discovery. For Plato, the collective unity of all these lines into a set, or a name, would be equally trivial, as he makes clear in the beginning of the Philebus 14d-e.[18]

---

[17] The proof of Proposition X.2 in the *Elements* makes use of Eudoxus' condition def.4, Book V, of the *Elements*. For this reason Knorr [K] believes that it could not be in use by Theodorus. Nevertheless an elementary, not involving Eudoxus' condition, proof of X.3 can be given with a simple argument, similar to the proof, of Proposition VI.1 of the *Elements,* suggested in Aristotle's *Topics* 158 b 29-35.

Both Knorr and Burnyeat translate the crucial sentence 147d7-8 in a collective, non-distributive way:
Knorr [Kn, 1975] (p.63): ‘Now this is what occurred to us: that, since we recognized the power to be unlimited in number, we might to collect them under a single name’;
Burnyeat [Bu, 1978]: ‘Well, since the powers seemed to be unlimited in number, it occurred to us to do something on these lines: to try to collect the powers under one term by which we could refer to them all.’

But they are both somewhat uneasy about the meaning of the sentence, and they settle for a similar interpretation. We will look at Burnyeat’s treatment:
‘The key sentence is 147d 7-e 1, which I render as follows:
*Since the dunameis were turning out to be unlimited in number, it occurred to us to attempt to collect them up into a single way of speaking [i.e., a formula or definition] of all these dunameis together.*
Theaetetus is recounting the thoughts suggested to himself and his companion by and during Theodorus' lesson, and the idea that there is an endless series of whole number squares (or sides of such squares) would hardly need to be prompted by a process as protracted as Theodorus' lesson. That there are an indefinite, perhaps infinite, number of squares with incommensurable sides, on the other hand, is precisely the hypothesis that would suggest itself as Theodorus proceeded from case to case proving more and yet more examples of incommensurability, perhaps by a method which could be endlessly reapplied. Therefore, it is likely that, in context, "all these dunameis refers to squares with incommensurable sides rather than to squares generally.’

We should point out that the fact that there are infinitely many powers is also obvious and not relate at all to Theodorus’ lesson, simce the incommensurability of the power for N=2, known to the Pythagoreans and considered known by Theodorus, implies the incommensurability of the infinitely many powers for $N=2^{2n+1}$, with n=1,2,….

[18] Cf. Section 14. b2 below

## 6. *Theaetetus* 148d4-5: 'peiro mimoumenos' the Theaetetean Division & Collection towards the Platonic Division & Collection in the *Sophistes* and the *Politicus*

The 'distributive' interpretation of the plural 'hai dunameis' in 147d7-8 opens the way for a meaningful understanding of But why was Plato interested for the Theaetetus mixture? It is not simply to eulogise Theaetetus; he has something more significant in mind.

### 6.1. The Plan to imitate Theaetetus' mathematical Division and Collection in philosophy

**(i) 145d6, 145e8-9**
Socrates, with a characteristic ironic understatement ('micron de ti aporo', 146d6) sets the basic question for the trilogy *Theaetetos, Sophistes, and Politicus,*
how is it possible to obtain knowledge of a Platonic Idea **('episteme').**

**(ii) 147c7-d1**
After some abortive attempts, Theaetetus realizes that
Socrates' question on the acquisition of knowedge of Platonic Ideas ('episteme')
is similar ('hoion'), relevant to
his (and his companion) treatment of the mathematical problem, conceived during Theodorus' lesson on quadratic incommensurabilities, and which was described in Section 5.**<** Every dunamis (Theaetetean mixture) possesses division and collection (*Theaetetus* 147d3-148b4)>.

**(iii) 148b5-8**
Theaetetus further realizes that, even though his answer for collecting a power into one is similar to the problem about episteme of Platonic Beings, still he cannot give the precise same answer for episteme of Platonic Beings, perhaps because the concept of dunamis involves not only commensurability and incommensurability (concepts that can be expressed in philosophical terms of finite and infinite) but also of taking **squares,** so that the geometric concept of dunamis has no clear analogue in the realm of Platonic Beings.
Thus one might consider, as in the Philebus later, that commensurability and incommensarability are really mathematical instances of the more comprehensive philosophical principles of finite and infinite, but still what sense is one to make in the more comprehensive philosophical context of **squaring**.[19]

**(iv) 148c9-d2**
*Socrates***.***Then you must have confidence in yourself (***tharrei** *peri sauto),*
*[148d] and try earnestly (***prothumetheti***) in every way to gain an understanding of the nature of knowledge as well as of other things.*
*Theaetetus***.** *If it is a question of earnestness, Socrates, the truth will come to light.*
**148d4-7**

[19] The philosophical analogue of squaring is the double Measurement in the *Politicus* 283a-287b, cf. Section 12, below.

*Socrates***.** *Well then—for you pointed out the way admirably just now—***try to imitate (peiro mimoumenos)** *your answer about the powers (dunameis), and just as you embraced (perielabes) them in one kind (heni eidei), though they were many (pollas), try to designate the many forms of knowledge (tas pollas epistemas) by one Logos.*

Theaetetus suceeded in showing that each power, **divided** into an infinite multitude of parts by Theodorus, can be **collected** into One; Socrates is exhorting Theaetetus to **imitate** ('peiro mimoumenos') this process of **Collection** to the knowledge of Platonic Ideas. This **imitation** is realized in the two remaining dialogues of the trilogy; in the *Sophistes* and *Politicus* the method of obtaining knowledge of the platonic Ideas is the method of Division and Collection ('diairesis kai sunagoge'), developed by the Eleatic Stranger in the presence of Theaetetus and his companion. This brings closer the description of Theaetetus mathematical achievement of **collecting the divided 'dunameis'** with the description, as **Collection and Division,** of the philosophical method of obtaining knowledge ('episteme') of platonic Ideas.

### **6.2**. The imitation realized with the method of Division and Collection in the dialogues *Sophistes* and *Politicus*

The two last dialogues, *Sophistes* and *Politicus*, of the Platonic trilogy on the question of knowledge of the Platonic Ideas, are dedicated to the method of Division and Collection by means of which humans are able to obtain complete knowledge of the intelligible Platonic Ideas. Thus
Theaetetus succeeded in collecting (sullabein, perilabein) the infinite multitude of the parts, into which a power had been divided, into One; and in the two succeeding dialogues Platonic Ideas become known by the method of Division ands Colletion.So the imitation of Mathematics to philosophy into which Socrates urged Theaetetus has been realized.
Since philosophical Collection in the *Sophistes* and the *Politicus* is an imitation of the geometrical Colletion in the Theaetetus, we will attempt to understand the double meaning of Collection by studying the method in these two dialogues. This will be done in Sections 9,10,11, and 12.

### 6.3. A fourth argument in favor of the distributive nature of the plural in the *Theaetetus* 147d7-8

In the *Sophistes* 247d-248d Plato defines the Platonic Idea as a 'dunamis'. The meaning of this definition will be explained briefly below in Section 14.4 (4). According to this definition the fundamental property of a Platonic Idea is that, in interaction, active or passive, with an opposite element, forming an indefinite dyad, it has the 'power', in fact the 'power to communicate' with each other (252d3, 253a8, 253e1, 254b8, 254c5), in fact the 'power for equalisation with each other' (257b).

Irrespetive of the precise meaning of the definition, Plato imitates here the term used in the *Theaetetus* 147d-148d for the mathematical 'dunamis'. In the *Sophistes* the Division and Collection is a method for the acquisition of knowledge of a Platonic Idea, thus of a 'dunamis'; it is certainly the case that every dunamis separately is Collected into One. This suggests again that in the *Theaetetus* the Collection of the 'dunameis' (who were collected into One because they appeared to be infinitely many) is to be understood distributively, for each power separately.

## 7. Some Plausibility arguments suggesting periodic anthyphairesis

The distributive interpretation of the *Theaetetus* 147d7-8, given in Section 6, and the resulting discussion about the nature of the imitation (*Theaetetus* 148d4-5) suggest some plausibility arguments that lead to periodic anthyphairesis.

### 7.1. First plausibility argument: Theodorus' and Theaetetus' method of proof of incommensurability of the powers is anthyphairetic

The fact that Theodorus' proofs of incommensurability of some powers resulted in the division of each of the powers **into an infinite multitude of parts** (according to the distributive interpretation of Section 5) is not consistent with any other conceivable method save an anthyphairetic one, eventually employing Proposition X.2 of the *Elements*. And the fact that Theaetetus continued from there on, collecting all these infinitely many pieces into one, indicates that Theaetetus' contribution was also anthyphairetic.

### 7.2. Second plausibility argument: Theaetetus' Collection leads to a complete knowledge of every 'dunamis'

The method of Theodorus is Division, while that of Theaetetius Collection (described by the ideas 'sullabein' and 'perilabein'); thus put together they deal with may be described as **Division and Collection**.
On the other hand Socrates urges Theaetetus to try to imitate the mathematical Division and Collection for the philosophical problem of the acquisition of knowledge of the Platonic Ideas; and, in later in the trilogy, in the *Sophistes* and *Politicus*, the method by which klowledge of the Platonic Ideas is acquired is described as **Division and Collection**.
Thus the philosophic, Platonic method of Division and Collection is an imitation of Theaetetus' (Division and) Collection.

The fact that the reason for Socrates interest in Theaetetus' mathematical discovery is the quest for **a method of acquiring knowledge** of the Platonic Beings indicates that Theaetetus' mathematical discovery about powers was not simply about a proof of their incommensurability but rather a method of acquiring knowledge of these powers. This again points to periodic anthyphairesis, since by periodic anthyphairesis of a power a complete knowledge of the power is indeed obtained, as indicated by the following theorem of modern mathematics.

### 7.3. The theorem on the continued fraction expansion of quadratic irationals

The following important theorem of modern mathematics ($17^{th}$ and $18^{th}$ century) shows that the **continued fraction expansion of a quadratic irrational** provides complete knowledge of the quadratic irrational.

**Theorem[20] (Fermat, Brouncker, Euler[21], Lagrange[22], Legendre, Galois[23])**.

---

20 cf. Hardy & Wright [HW, 1938], Weil [W, 1984], Fowler [Fow, 1986].
21 [E,1765]
22 [Lag, 1769], [Lag, 1771]
23 [Ga, 1829]

If N is a non-square natural number, then the continued fraction of √N is periodic, and in fact palindromically periodic, namely there are numbers n and m, $k_1, k_2, \ldots, k_{n-1}, (k_n)$, such that the sequence of successive quotients of the continued fraction of √N has the form

**[m, period ($k_1, k_2, \ldots, k_{n-2}, k_{n-1}, (k_n), k_{n-1}, k_{n-2}, \ldots, k_2, k_1$,2m)]**.

The modern theory of continued fraction of a real number is the equivalent analogue of the anthyphairesis of a ratio of, say, line segments. The above theorem is stated equivalently as follows:

**Proposition**. If two lines a and b are incommensurable in length only, namely form a Theaetetean mixture of incommensurability and commensurability, then the anthyphairesis of a to b is **periodic**, and in fact **palindromically periodic**. In fact, there are numbers n and m, $k_1, k_2, \ldots, k_{n-1}, (k_n)$, such that the sequence Anth(a,b) of successive quotients of the continued fraction of a to b has the form Anth(a,b)=[m, period ($k_1, k_2, \ldots, k_{n-2}, k_{n-1}, (k_n), k_{n-1}, k_{n-2}, \ldots, k_2, k_1$,2m)].

Thus the anthyphairesis of the Theatetean mixture **provides complete, full knowledge** of the Theaetetean .mixture

**7.4. Third plausibility argument: The cyclical-periodic nature of the method of Collection**

Theaetetus' 'collection into one', a mathematically non-familiar notion, initially described by 'sullabein eis hen' (147d8), is also further described by **'heni eidei perielabes'** (148d6), a term that has a cyclic and possibly periodic connotation.

On the other hand there are several descriptions of the method of Division and Collection in Platonic dialogues, specifically in the *Sophistes*, *Politicus*, *Phaedrus*, and *Philebus*. It is remarkable that in every one of these descriptions, Collection is described Before embarking on a detailed study of the method we will stop to note that all the descriptions of the (Division and) Collection method imply a by an expression that denotes that it has a cyclical-periodic nature.[24]

---

[24]Here is a list of such descriptions, indicating the term, in each description showing cyclical nature-periodicity.
*Sophistes* 235a10-c6 ('[ton thera] **perieilephamen en amphiblestrikoi tini'**);
*Sophistes* 253d1-e3 ('pollas heteras allelon hupo mias exothen **periechomenas'**);
*Sophistes* 264d12-265a2 ('schizontes dichei to protethen genos, poreuesthai kata toupi dexia aei meros tou tmethentos, echomenoi tes tou sophistou koinonias, heos an autou ta koina panta **perielontes**, ten oikeian lipontes phusin epideixomen');
*Politicus* 283b1-c2 ('**perielthomen en kukloi** pampolla diorizomenoi');
*Politicus* 285a4-b6 ('tas de au pantodapas anomoiotetas…sumpanta ta oikeia **entos mias homoiotetos herxas** genous tinos ousiai **peribaletai'**);
*Politicus* 286d6-287a6 ('logon meke kai **tas en kukloi periodous'**);
*Philebus* 15d4-8 ('tauton hen kai polla hupo logon gignomena **peritrechein** pantei kath' hekaston ton legomenon aei').
In *Phaedrus* 264e-266c the Division and Collection of the right Love is periodic by its opposition to the explicit non-periodicity of the left-sinister Love, indicated bythe statement 'to ep' aristera temnomenos meros, palin touto temnon **ouk epaneken prin en autois'**. (The sentence usually receives a different reading, in which periodicity is lost, e.g. Hackforth [Ha, 1945]).

# 8. The Logos Criterion for periodic anthyphairesis

## 8.1. The Theaetetean theory of ratios of magnitudes

Aristotle, in the passage 158-9 of the *Topics*, important for the history and reconstruction of Greek mathematics, informs us that before the theory of ratios of magnitudes that dominated in Euclidean geometry, and which is Eudoxus' theory majestically presented in Book V of the *Elements*, there was a previous theory based on the

**Definition**. If a,b are homogeneous magnitudes and b,d are homogeneous magnitudes then
a/b=c/d if and only if Anth(a,b)=Anth(c,d)
(where Anth(a,b) denotes the sequence of successive quotients of the anthphairesis of a to b).

A theory of ratios based on this pre-Eudoxian definition would need Eudoxus' condition def. 4 in Book V. of the *Elements*, and so would be historically untenable. (Knorr [Kn, 1975], appendix B, reconstructs such a theory (and as carefully noted by Knorr, the crucial Theorem 6, p. 338, the analogue of the fundamental Eudoxian Proposition V.8 of the *Elements*, makes essential use of the Eudoxus condition (V, def. 4). There is however evidence, chiefly among some that this pre-Eudoxian theory applied only to a limited class of pairs of magnitudes (e.g. the rational lines wrt an assumed line in the sense of definition 3 of Book X of the *Elements*—hence the term 'alogos' for those ratios whose square is not commensurable), which was quite adequate for the purposes of quadratic incommensurabilities, and for which there is no need to an appeal to Eudoxus condition. One anonymous coment in *Scholia in Eucliden* V.22 suggests that there was such an earlier restricted quadratic theory of magnitudes, and another (*Scholia in Eucliden* X.2) criticizes the general Eudoxian theory, clearly from a Platonic point of view and in favor of the earlier quadratic theory of ratios of magnitudes.
It is generally accepted that the theory of ratios based on anthyphairesis is due to Theaetetus.

## 8.2. The Logos Criterion for periodic anthyphairesis

A natural consequence of Theaetetus' definition of equality of ratios for magnitudes is the following

**Proposition (Logos Criterion for the periodicity of an anthyphairesis)**.
Let a,b be two line segments, with a>b, and anthyphairesis given as follows:

$a = I_0 b + e_1,$ with $b>e_1,$
$b = I_1 e_1 + e_2,$ with $e_1>e_2,$
…
$e_{n-1} = I_n e_n + e_{n+1},$ with $e_n > e_{n+1},$

$e_n = I_{n+1} e_{n+1} + e_{n+2},$ with $e_{n+1}>e_{n+2},$
…
$e_{m-1} = I_m e_m + e_{m+1},$ with $e_m > e_{m+1},$

$e_m = I_{m+1} e_{m+1} + e_{m+2},$ with $e_{m+1}>e_{m+2},$
…
…….

and assume that there are indices n<m such that

$e_n/e_{n+1} = e_m/e_{m+1}$ (‘Logos Criterion’).

Then the anthyphairesis of a to b is eventually periodic, in fact
$Anth(a, b) = [I_0, I_1, …, In, period(I_{n+1}, I_{n+2}, …, I_m)]$.

**Proof**. By the Theaetetean definition of analogy, $Anth(e_n, e_{n+1}) = Anth(e_m, e_{m+1})$.
But Anth (a,b)=
$[I_0, I_1, …, I_n, Anth(e_n, e_{n+1})]=$
$[I_0, I_1, …, I_n, I_{n+1}, I_{n+2}, …, I_m, \mathbf{Anth(e_m, e_{m+1})}]=$
$[I_0, I_1, …, I_n, I_{n+1}, I_{n+2}, …, I_m, \mathbf{Anth(e_n, e_{n+1})}]=$
$[I_0, I_1, …, I_n, I_{n+1}, I_{n+2}, …, I_m, \mathbf{I_{n+1}, I_{n+2}, …, I_m, Anth(e_m, e_{m+1})}]=$
…=
$Anth(a, b) = [I_0, I_1, …, I_n, \mathbf{period(I_{n+1}, I_{n+2}, …, I_m)}]$.

**8.3.** An **abbreviated representation of the Logos Criterion**

We note that the Logos criterion is stated solely in terms of the remainders, and no use is made in its statement of the quotients. The Logos criterion can be adequately represented in the following abbreviated form:

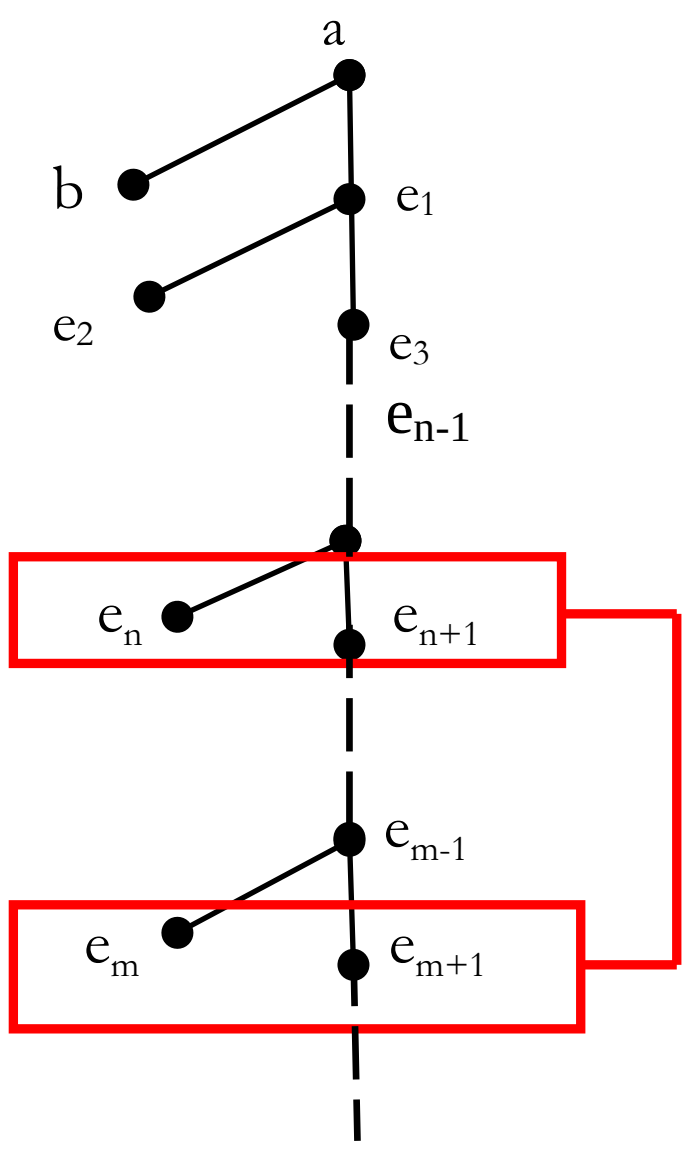

**8.4. Example. The anthyphairesis of a to b with $a^2=19b^2$**

We exhibit the logos criterion with the first power that Theodorus did not prove in his lesson to Theaetetus and his companion.

**(i) The anthyphairesis of a to b**

$a=4b+e_1$, με $a_1<b$ (whence $e_1=a-4b$), (and $b=2e_1+e_2$, $e_2<e_1$ (whence $e_2=9b-2a$));
$e_1=e_2+e_3$, $e_3<e_2$ (whence $e_3=3a-13b$),(and $e_2=3e_3+e_4$, $e_4<e_3$ (whenc $e_4=48b-11a$)),
$e_3=e_4+e_5$, $e_5<e_4$ (whence $e_5=14a-61b$), (and $e_4=2e_5+e_6$, $e_6<e_5$ (whence $e_6=170b-39a$)),
$e_5=8e_6+e_7$, $e_7<e_6$ (whence $e_7=326a-1421b$).

**(ii) Verification of the Logos Criterion**
**$b/e_1=e_6/e_7$.**
(we simply verify that $b.e_7=b.(326a-1421b)=e_1.e_6=(a-4b).(170b-39a)$).

The Logos Criterion results in the full knowledge of the ratio a/b:
Anth (a,b)= [4, period(2,1, 3,1,2,8)].
By Proposition X.2, of course it is concluded that a,b are incommensurable.

**(iii) Abbreviated representation of the Logos Criterion**

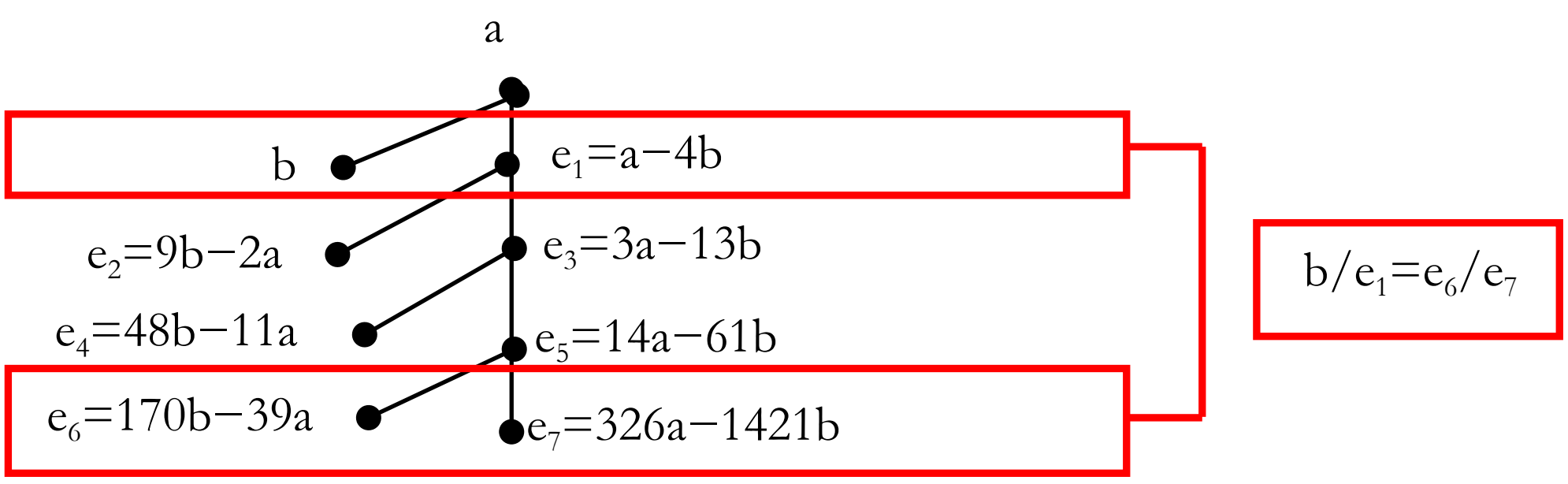

# 9. The Division and Collection of the Platonic Idea ‘the Angler’ in the *Sophistes* 218b-221c

## 9.1. The Division of the ‘Angler’

The method of Division and Collection, also called ‘Name and Logos’ (cf. *Theaetetus* 201e2-202b5, *Sophistes* 218c1-5, 221a7-b2, 268c5-d5), is exemplified at the start of the *Sophistes* 218b-221c by the definition of the Being Angler. In the scheme below, we reproduce the binary division process which leads to the Angler.

Table 9.1. The Division of the Angler (*Sophistes* 218b-221c)

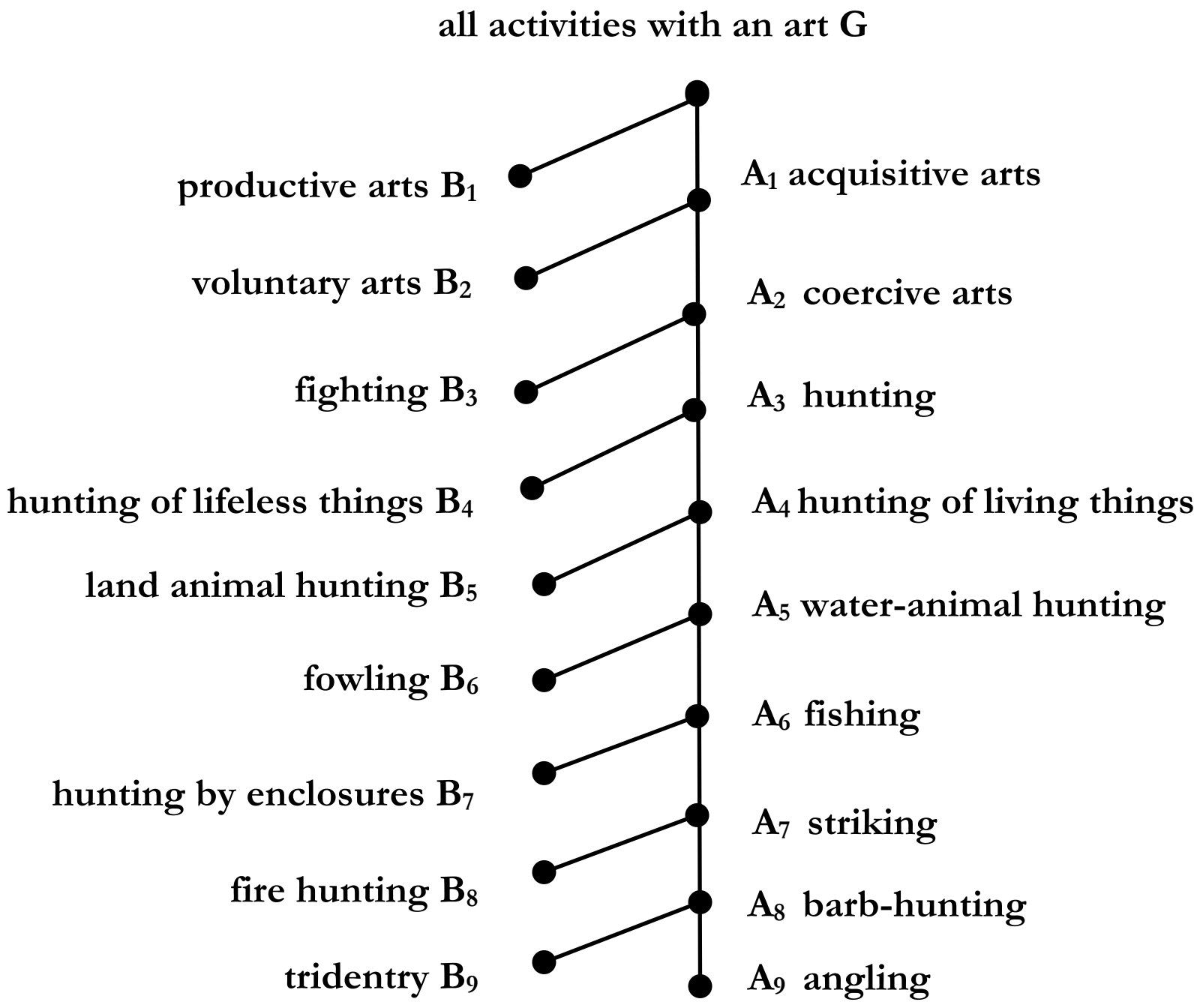

## *9.2. Collection-Logos of the Angler*

The description of the Logos-Collection of the Angler is contained in the *Sophistes* 220e2-221c3 passage (which we have divided into two parts [A] and [B] for the sake of convenience):
[A] ‘Stranger: Then of striking which belongs to barb-hunting [A8], that part which proceeds **downward from above** (‘anothen eis to kato’),
is called, because tridents are chiefly used in it, tridentry [B9],
I suppose….
Stranger: The kind that is characterized by the opposite sort of blow, which is practised with a hook and strikes,... and proceeds
**from below upwards** (‘katothen eis tounantion ano’),
being pulled up by twigs and rods.
By what name, Theaetetus, shall we say this ought to be called?
Theaetetus: I think our search is now ended and we have found the very thing we set before us a while ago as necessary to find.
Stranger: Now, then, you and I are
only agreed about the name of angling,[A9] ’ (220e2-221b1)

[B] ‘but **we have acquired also a satisfactory ‘Logos’** of the thing itself.
**For (‘gar’)**
of art as a whole, half was acquisitive,
and of the acquisitive, half was coercive,
and of the coercive, half was hunting,
and of hunting, half was animal hunting,
and of animal hunting, half was water hunting,
and of water hunting [A5]
**the whole part from below (‘to katothen tmema holon’)**
was fishing,[A6]
and of fishing, half was striking,
and of striking, half was barb-hunting,[A8]
and of this [A8]
the part in which the blow is pulled
from below upwards (‘peri ten katothen ano’) was angling.[A9]’[25] (221b1-c3).
In [A], the opposing relation of Tridentry to Angling is carefully explained: all Fishing with a hook is divided into
Tridentry (=Fishing with a trident), which is described as
Fishing with a hook with an art that proceeds **from above downwards**, and
Angling (=Fishing with a rod), which is described as
Fishing with a hook with an art that proceeds **from below upwards**.
We have pushed Division all the way to the Angling; thus, we have certainly found ‘the name’ of Angling.

But now, in [B], it is claimed that ‘the Logos’ of the Angling has been found, too.
The justification—the proof—that we have indeed found the Logos, too, is contained in the remainder of [B], since this remaining part of [B] starts with a ‘for’ (‘gar’), and this justification can be seen to consist of:
(i) an accurate recounting of all the division-steps, abbreviated in the sense that of the two species into which each genus is divided, only the one that contains the Angler is mentioned, while the species opposite to it is omitted;
(ii) a reminder that the last species, the angling, is characterised as the part of its genus that proceeds ‘from below upwards’; and,
(iii) the ONLY new information (since (i) and (ii) are repetitions of things already contained in the Division and in [A]), which concerns the species of fishing, three steps before angling, and which informs us for the first time that this species is ‘the whole part from below’ of its genus.
Since this is an abbreviated account, there is no explicit information on the species opposite to ‘fishing’—namely ‘fowling’, but since ‘fishing’ was described not simply as ‘the part from below’ of its genus, but emphatically as ‘the whole part from below’, it follows that the opposite species—‘fowling’—must be characterised as ‘the (whole) part from above’ of the same genus. In fact, there can be no other justification for the presence of the term ‘whole’ in the description of ‘fishing’ with a view to arguing that we have obtained ‘Logos’, except to indicate and imply this description for its opposite species, ‘fowling’.
We recall that the part of [B] from the word ‘for’ (‘gar’) is explicitly a justification of the claim that we have succeeded in finding the ‘Logos’ of the Angling. We may then ask: what is the 'Logos' of the Angler that reasonably results from such a justification? There can really be only one answer: the 'Logos' we are looking for is the equality of the 'philosophic ratio' of Tridentry to Angling, namely the equality of the ratio 'of from above downwards to from below upwards', to the ratio 'of Fowling to Fishing'.

[25] Based on Plato, *Theaetetus and Sophist*, translated by H.N. Fowler, Loeb Classical Library, Cambridge, Mass., 1921.

Thus [9b]

**fowling B6/fishingA6=**
**tridentry B9/Angling A9=**
**from above downwards/from below upwards.**

Since the species Tridentry and Angling forms a pair of opposite species, and the species Fowling and Fishing form another pair of opposite species in the Division Scheme for the Angler, the resulting 'Logos' bears a most uncanny similarity to the Logos Criterion for the periodicity of the anthyphairesis of geometric magnitudes, and of geometric powers in particular.

### 9.3. The Division and Collection of the Angler

The Division and Collection of the Angler thus takes the following form:

Table 9.3. Division and Collection of the Angler (*Sophistes* 218b-221c)

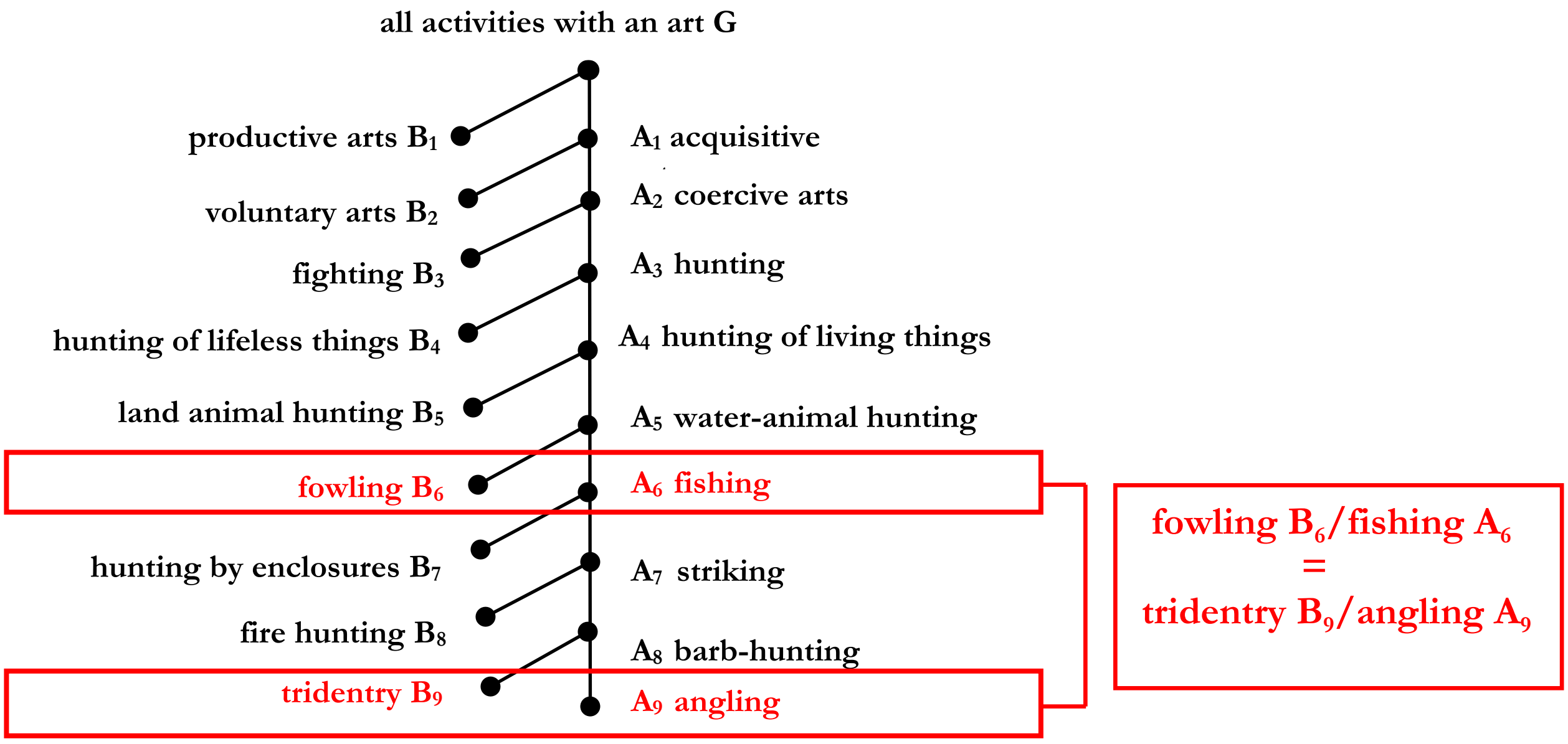

Thus the Division and Collection of the Angler consists of the Division, described in a, which is analogous to the abbreviated anthyphairetic model, as given in 1a, and of the Logos, described in b, analogous to the Logos Criterion for periodicity for geometric anthyphairesis established in Section 5.

## 10. The Division and Collection of the Platonic Idea ‘the Sophist’ in the *Sophistes* 234e-236d & 264b-268d

It now appears that the Division and Collection of a Platonic Idea—and the Angler is certainly a lowly paradigm of a Platonic Idea—is very much like the anthyphairetic Division and the Logos Criterion of a geometric ‘power’. When Socrates expressed, in the *Theaetetus* 145c7-148e5, his exhortation to imitate the geometric situation, it seems that he meant a much closer imitation that anybody has thus far suspected! But before proceeding on to wide-ranging conclusions, it would be prudent to examine whether the Division and Collection of the Sophist, in the *Sophistes* 264b-268d, has the same kind of structure and, in particular, whether there is a similar type of ‘Logos’. So we shall take a look at the Division and Collection of the Sophist.

### 10.1. The Division of the Sophist

The Division for the Sophist follows the same pattern as the Division of the Angler, starting with a Genus, in this case ‘all the productive arts’, proceeding by binary division of each Genus to two species, where the next Genus is that species of the previous step in the Division that contains the entity to be defined, in this case the Sophist, and ending with the division-step that produces the Sophist as a species. The whole division scheme is as follows:

Table 10.1. Division of the Sophist (*Sophistes* 264b-268d)

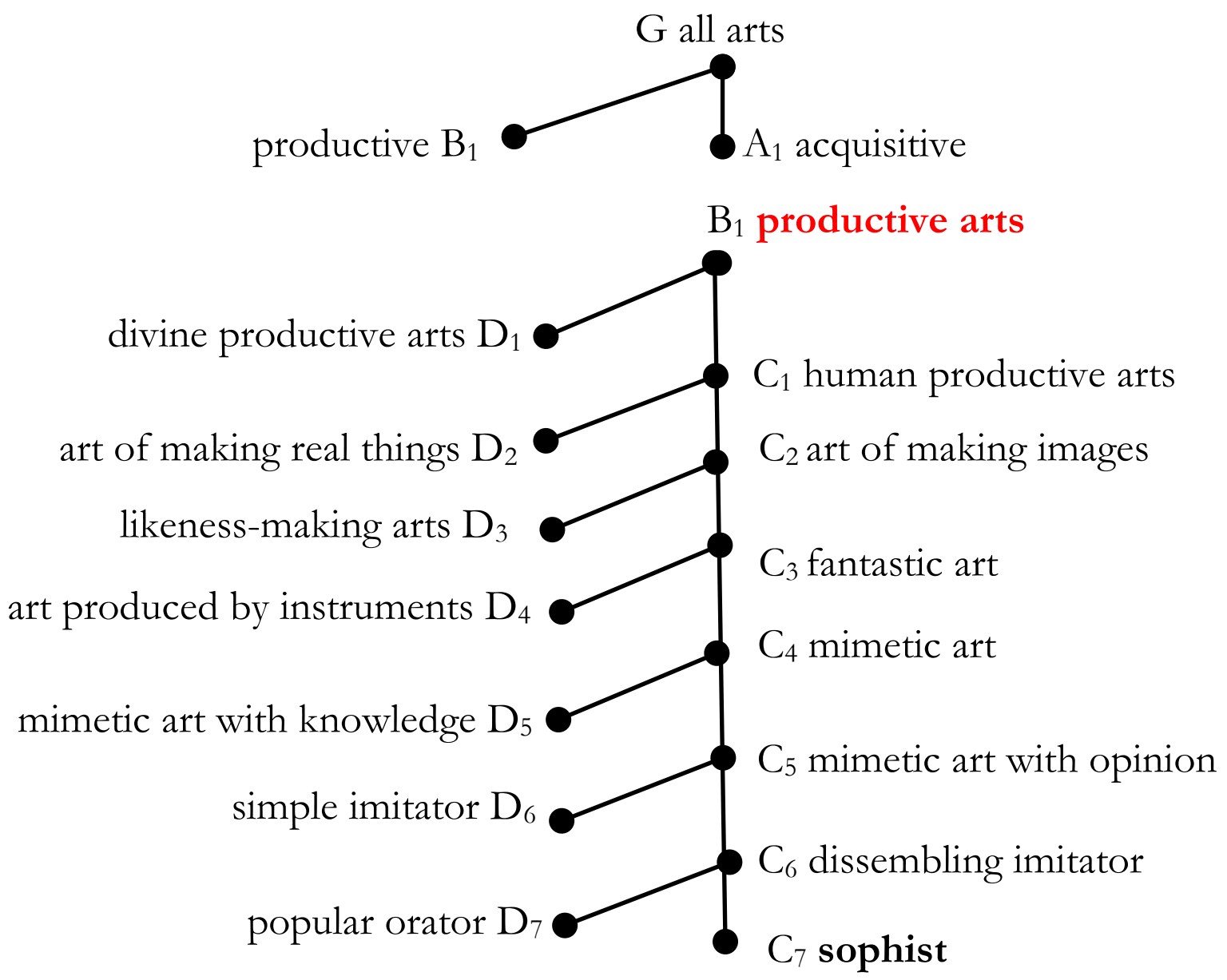

### 10.2. The fundamental analogy in the Divided Line of the *Politeia* 509d-510b and the Logos criterion-Collection of the Sophist

We will now prepare the ground for the Logos-Collection of the Sophist. The fundamental analogy in the Divided Line of the *Politeia* 509d-510b plays a central role in the Logos Criterion of the Sophist. Here is the passage:

<table>
<tr><td colspan="3">‘Conceive then, said I, as we were saying, that there were two entities, and that</td></tr>
<tr><td colspan="2"></td><td>one of them is sovereign over<br>the intelligible order and region</td></tr>
<tr><td colspan="2">and the other over<br>the world of the eyeball</td><td></td></tr>
<tr><td colspan="3">…You surely apprehend the two types,</td></tr>
<tr><td colspan="2">the visible (‘horaton’)</td><td></td></tr>
<tr><td colspan="2"></td><td>and the intelligible (‘noeton’)</td></tr>
<tr><td colspan="3">…Represent them then, as it were, by a line divided into two unequal sections<br>and cut each section again in the same ratio (‘ana ton auton logon’)</td></tr>
<tr><td colspan="2">(the section, that is, of the visible</td><td></td></tr>
<tr><td colspan="2"></td><td>and that of the intelligible order),</td></tr>
<tr><td colspan="3">and then as an expression of the ratio of their comparative</td></tr>
<tr><td colspan="2"></td><td>clearness</td></tr>
<tr><td colspan="2">and obscurity</td><td></td></tr>
<tr><td colspan="3">you will have, as</td></tr>
<tr><td>one of<br>the sections of<br>the visible world,<br>images (‘eikones’).<br>By images<br>I mean, first,<br>shadows, and their<br>reflections<br>in water and on<br>surfaces of dense,<br>smooth and bright<br>texture, and<br>everything<br>of that kind…</td><td></td><td></td></tr>
<tr><td></td><td>As the second<br>section assume<br>that of which this is<br>a likeness<br>or an image<br>(‘ho touto eoike’),<br>that is, the animals<br>about us and all<br>plants and the<br>whole class of<br>objects<br>made by man…</td><td></td></tr>
<tr><td colspan="3">Would you be willing to say, said I, that the division in respect of</td></tr>
<tr><td colspan="2"></td><td>truth</td></tr>
<tr><td colspan="2">or the opposite</td><td></td></tr>
</table>

| | |
|---|---|
| | (‘aletheia te |
| kai me’) | |
| is expressed by the proportion:<br>as is | |
| the opinable (‘doxaston’) | |
| | to the knowable (‘gnoston’) |

so is

| | | |
|---|---|---|
| the likeness<br>(‘to homoiothen’) | | |
| | to that of which<br>it is a likeness<br>(‘to ho omoiothe’)? | |
| I certainly would.’[26] | | |

This analogy in the *Politeia* Divided Line is rendered as follows:

| | |
|---|---|
| Let L be a line (segment),<br>and divide the line L into two unequal sections,<br>say | |
| | A |
| and B, | |

with

| | |
|---|---|
| | A<br>representing the intelligible domain, |
| and B<br>the visible, or rather the sensible, domain. | |

(Incidentally, the construction of the division of a line segment into a given ratio is contained in Proposition VI.10 of the *Elements*).

| | | | |
|---|---|---|---|
| | | Then divide<br>section A<br>into two sections, say | |
| | | | C |
| | | and D, | |
| and divide<br>section B<br>into two sections, say | | | |
| | E | | |
| and F, | | | |
| in such a way that B/A=D/C=F/E. Further, | | | |
| F<br>represents<br>**the images** | | | |
| in the sensible domain B, | | | |

[26] Based on Plato, *The Republic*, translated by P. Shorey, Loeb Classical Library, Cambridge, Mass.,1935.

<table>
<tr><td></td><td>and E<br>the entities</td><td></td></tr>
<tr><td colspan="2">in the sensible domain B<br>to which</td><td></td></tr>
<tr><td>these are images,</td><td></td><td></td></tr>
<tr><td></td><td>i.e.<br>the real entities</td><td></td></tr>
<tr><td colspan="2">in the sensible domain B.</td><td></td></tr>
<tr><td colspan="2"></td><td>Also the intelligible domain A<br>is identified with<br>the domain of knowledge,</td></tr>
<tr><td colspan="2">and the sensible domain B<br>is identified with<br>the domain of opinion.</td><td></td></tr>
<tr><td colspan="3">Hence the following proportion holds:<br>The ratio of the opinable B to the knowable A<br>is equal to<br>the ratio of the images F to the real entities E.</td></tr>
</table>

This, then, is the fundamental analogy of the Divided Line in the *Politeia* 509d-510b:

**the opinable/the knowable = the likeness/that of which it is a likeness**.

Therefore the Logos Criterion is now immediately recognised:

**the art of making real things D2/ the art of making images C2 =**
**mimetic art with knowledge D5/mimetic art with opinion C5**.

**10.3. The next two Logoi of the Division and Collection of the Sophist analogous to the geometrical ratios of periodic anthyphairesis.**

The definition of the Sophist presents an additional feature that ties it even closer to the mathematical model: there are two Logoi after the Logos Criterion, and if the mathematical anthyphairetic model is indeed followed, we would expect to have two further equalities of Logoi: namely the logos of the third step should be equal to the ratio of the sixth division-step, and the logos of the fourth division-step should be equal to the ratio of the seventh and final division-step**.**

(c1) The equality of the third and sixth ratios of the Division of the Sophist. In the third division-step, the Genus ‘humanly produced images’ is divided into two species, ‘Stranger: I see the likeness-making art (‘eikastiken’) as one part of [the image-making art]. This is met with, as a rule, whenever anyone produces the imitation by following the proportions of the original in length, breadth, and depth, and giving, besides, the appropriate colors to each part.’ (235d6-e2);
‘Stranger: Now then, what shall we call that which appears, because it is seen from an unfavorable position, to be like the beautiful, but which would not even be likely to resemble that which it claims to be like, if a person were able to see such large works adequately? Shall we not call it, since it appears, but is not like, an appearance?

Theaetetus: Certainly.’ (236b);
‘Stranger: And to the art which produces appearance, but not likeness (‘phantasma all’ ouk eikona’), the most correct name we could give would be "fantastic art," (‘phantastiken’) would it not?
Theaetetus: By all means.
Stranger: These, then, are the two forms of the image-making art (‘eidolopoiikes’) that I meant, the likeness-making (‘eikastikes’) and the fantastic (‘phantastikes’).’ (236c3-7);
‘Stranger: We must remember that there were to be two parts of the image-making class (‘eidolourgikes’), the likeness-making (‘eikastikon’) and the fantastic (‘phantastikon’)’[19] (266d8-9).

In the sixth division-step, described in 267e7-268a8, the Genus opinionated-imitator is divided into the two species, ‘the simple’ (‘haploun’) and ‘the dissembling’ (‘eironikon’) opinionated-imitator, the simple imitators being those who
‘think they know that about which they have only opinion’,
while the ‘dissembing’ are those who
‘because of their experience in the rough and tumble of arguments, strongly suspect and fear that they are ignorant of the things which they pretend before the public to know’.

Thus the simple imitators do not distort their opinion, but rather express a likeness of their opinion, while the dissemblers distort and disguise their opinion behind a false appearance. Therefore we have

**Proposition**. The ratio of the third step in the Division of the Sophist, namely **the ratio of the likeness-making arts to the phantastic arts**, is equal to the ratio of the sixth step in the Division of the Sophist, namely **the ratio of the simple imitator to the dissembling imitator**.

(c2) The equality of the fourth and seventh ratios. In the fourth division-step, described in 267a1-b3, the Genus, fantastic arts, is divided into two species, as follows:
‘Stranger: Let us, then, again bisect the fantastic art.
Theaetetus: How?
Stranger: One kind is that produced by instruments, the other that in which the producer of the appearance offers himself as the instrument.
Theaetetus: What do you mean?
Stranger: When anyone, by employing his own person as his instrument, makes his own figure or voice seem similar to yours, that kind of fantastic art is called mimetic.’[19]
Thus, in mimetic art the instrument of imitation is the imitator himself, while in the nameless opposite art the instrument of imitation is other than the imitator.
In the seventh division-step, described in 268a9-c4, the dissembler is divided into the demagogue and the sophist, the demagogue being
‘one who can dissemble in long speeches in public before a multitude’,
while the sophist is someone
‘who does it in private in short speeches and forces the person who converses with him to contradict himself.’

Thus, he who listens to a dissembler is deceived, if that dissembler is a demagogue, and contradicted not by himself but by another instrument of deceit (namely the demagogue himself), while if the dissembler is a sophist, the listener is forced by the sophist to himself become the instrument of deceit.
Thus we have:

**Proposition**. The ratio of the fourth step in the Division of the Sophist, namely **the ratio of the imitator who uses other instruments of imitation to the imitator who is himself the instrument of imitation**
is equal to the ratio of the seventh step in the Division of the Sophist, namely **the ratio of the demagogue, whose listener is contradicted and deceived not by himself but by another, to the sophist, whose listener is forced to be contradicted and deceived by himself.**

Propositions (c1) and (c2) provide powerful additional evidence in favor of the anthyphairetic interpretation of Division and Collection.

**10.4. The Division and Collection of the Sophist. The complete Division and Collection of the Sophist can thus be summarised in the following**

Table 10.4. Division and Collection of the Sophist (*Sophistes* 264b-268d)

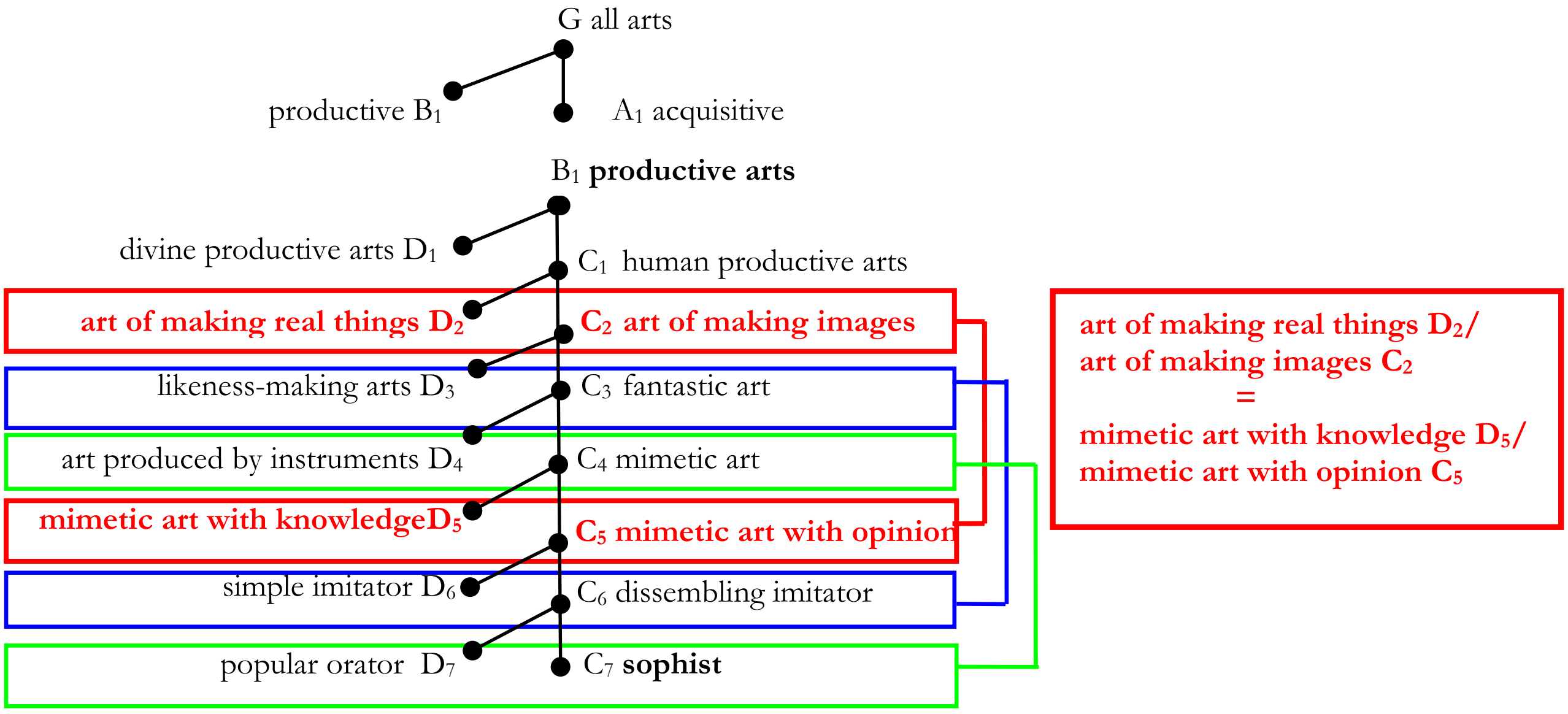

Thus, we have full confirmation for the anthyphairetic meaning of Logos that we proposed in the Division and Collection of the Angler. This time, the Logos is not of the simple-minded, obvious type that appeared in the definition of the Angler (from above downwards/from below upwards), but of a more sophisticated and philosophical kind, playing a central role in the dialectics of the *Politeia*[27].

[27] The correlation of the *Politeia* ratio in the divided line with the Division and Collection of the Sophist in the *Sophistes* is hinted at by Proclus in *eis Politeian* 1,290,7-10.

### 10.5. Concluding note for Setions 9 and 10.

We have argued that each of the two examples of Division and Collection contained in the *Sophistes*, that of the Angler and of the Sophist, has the structure of a philosophic analogue of periodic anthyphairesis.
Although the arguments leading to this interpretation are strong, one might remark that there is no indication to be found anywhere in the *Sophistes* that the Division is anything but finite. On the contrary, the language employed suggests that the Division ends when we arrive at an **indivisible** species. The language appears to contradict our claim about infinity and periodicity of the Platonic Division. Nevertheless we shall see in the *Politicus* Plato definitely describes the Division and Collection as periodic (Section 11c). However the question on indivisibility needs be explained. This will be done in Section 14.

## 11. The Division and Collection of the Platonic Idea ‘the Statesman’ in the *Politicus*

The third and final dialogue in the trilogy on Plato’s theory of knowledge is the *Politicus*. We will outline the contents of this dialogue, proceeding, not in the order presented by Plato but, in logical order.
The art of statesmanship is defined by the method ‘name and logos’, equivalently ‘division and collection’, involving steps with binary division and the logos criterion, in agreement with the method employed in the *Sophistes*.
We have seen in Sections 9 and 10 that in the *Sophistes* a Platonic Idea, like the Angler or the Sophist, becomes known to us humans by the definition of an Art, the Art of Angling or the Art of Sophistry, respectively, employing the method of division and collection, equivalently of ‘name and logos’, which we interpreted as a philosophic analogue of periodic anthyphairesis. It follows, from 11a below, that the Idea Statesman in the *Politicus* becomes known to us by the definition of the art of Statesmanship, employing the method of division and collection, equivalently of ‘name and logos’, a method already interpreted as the philosophic analogue of the eventually periodic anthyphairesis.

### 11.1. The first ten steps of the Division and Collection (258b-267c)

The aim of the *Politicus* is to obtain the knowledge (episteme) of the Platonic Idea of the Statesman, by the method of Division and Collection, the method employed in the *Sophistes*. In the first part of the dialogue, 258b-267c, the first ten steps of the Division are obtained. The initial genus G of all arts is successively divided in the binary manner familiar from the *Sophistes*. Thus

G divided into (K1, L1),
L1 divided into (K2, L2),
…,

L9 divided into (K10, L10).

Here is the Figure summarising these ten initial steps of the Division.

| **[G] all arts (258b-259d)** | |
|---|---|
| **[K1]** practical | **[L1]** scientific (‘gnostike’) |

| **[L1]** scientific (259d-260c) | |
|---|---|
| **[K2]**appraise, determine, **(**‘krinein’) | **[L2]** command, |
| | |

| **[L2]** command |
|---|

| (260c-261b) | |
|---|---|
| **[K3]** someone else’s commands | **[L3]** one’s own commands |
| | |

| **[L3]** one’s own commands (261b-d) | |
|---|---|
| **[K4]** lifeless herd | **[L4]** herd of living beings (=animals) |
| | |

| **[L4]** herd of animals (261d-e) | |
|---|---|
| **[K5]** herd of a **single** animal (‘monotrophia’) | **[L5]** herd of many animals |
| | |

| **[L5]** herd of many animals (264b-d) | |
|---|---|
| **[K6]** herd of aquatic animals | **[L6]**herd of land animals |
| | |

| **[L6]** herd of many land animals (264e-265b) | |
|---|---|
| **[K7]** herd of flying animals | **[L7]** herd of walking animals |
| | |

| **[L7]** herd of many walking land animals (265b-d) | |
|---|---|
| **[K8]** herd of horned animals | **[L8]** herd of hornless animals |
| | |

| **[L8]** herd of many hornless walking land animals (265d-e) | |
|---|---|
| **[K9]** herd of mixing animals | **[L9]** herd of non-mixing animals |
| | |

| **[L9]** herd of many non-mixing hornless walking land animals (265e-266d) |
|---|

| [**K10]** herd of four-footed (diameter's diameter) **pigs** | **[L10]** herd of two-footed (diameter) (human) |
|---|---|
| | |

Note that after L4 we are dealing with
**the arts of feeding an animal herd ('agelaiotrophike',** 261b-d**)**;
what is really being divided is the nature of the herd, narrowing it down, by successive binary division and exclusion, till we get to
**the arts of commanding a human herd (**L10, 265e-266d).

**11.2. The argument for the need to continue the Division** (267c-268d and 274e-277c)

The division steps described in 11a is not the end of our quest for the definition of the Statesman, as one might think. This is so because L10 consists of
**the arts of commanding a human herd**
but we are interested to define
the statesmanship,
which is
**the optimal art of commanding a human herd**.
This means that there is still ahead of us a second process of successive binary divisions, in which the herd is fixed as a human herd, and what is being now narrowed down, by binary division and exclusion, is
**the nature of commanding this human herd ('agelaiokomike',** 275e**)**,
till we get to
**the optimal commanding**.

**11.3. The myth of palindromic periodicity with Kronos era and Zeus era** (268d-274e)

According to this myth 268d-274e, there are two eras in the world, the Kronos era and the Zeus era.
The world, to be thought of as spherical entity, rotates in a certain direction during the Kronos era. After a long time, when all the possible stages of the Kronos era are exhausted ('to geinon hede pan aneloto genos' 'pasas tas geneseis apodedokuias'), it is necessary to reverse the direction ('palin anestrephen') and to go from the Kronos era to the Zeus era (272d6-e6).[28]
When the Zeus era is in danger of losing its direction and God perceives that it is in trouble, and that 'it might founder in the tempest of confusion and sink in the infinite sea of dissimilarity ('kedomenos hina me heimastheis hupo taraches dialetheis eis ton tes anomoiotetos apeiron onta ponton'), he … reversed ('strepsas') whatever had become unsound and unsettled in the previous period when the world was left to itself, set the world in order, restored it and made it immortal ('athanaton') and ageless (‘ageron’)’, thus returning (273c4-e4). In this way, we are back again to the Kronos period and completing a full period.[29]

[28] In a reconstruction of the proof of the theorem that every Theaetetean mixture has a palindromically periodic anthyphairesis, given in Section 13, there is a pigeonhole argument exactly at thye end of the semiperiod, namely in Plato;s terminology at the point of change from the Kronos era to the Zeus era. The reason given above for the change (aneloto, apodedokuias) appears to fit well with the pigeonhole argument.

The world moves forever repeating this period, alternating the Kronos era with the Zeus era.

The relation of the two eras is especially interesting. It is chematically described as follows: In the Zeus era (in which era we are presently) men come from non-existence a into existence b, then they grow c, then they get a beard d, then their beard becomes white e, then they grow old and are about to die f. If at this time the world happens to switch to the Kronos era, then the human will go from stage f to e then to d then to c then to b and then they pass into nonexistence a (270d-271a).

Summarising, according to the myth, the world is periodic, and the period consists of two half periods, the Kronos era followed by the Zeus era, in other words it is palindromically periodic..

**11.4. Every Division and Collection has palindromically periodic form, and is thus an infinite division. In particular the first ten steps of the Division and Collection of the Statesman correspond to the Kronos era and must be completed with the divisions of the Zeus era, in palindromic to the Kronos era manner, to be followed by the embracement of both in the Collection by means of Logos** (276a)

The passage 276a1-7 is crucial for our interpretation of the *Politicus*, because in this passage Plato makes perfectly clear that the palindromic periodicity of the myth must be applied for the Division and Collection of the Statesman, and since the Statesman is merely an example of an Art, and its Division and Collection is meant to show the general method, it certainly follows that every Division and Collection (including that of the Angler [Section 9], that of the Sophisth [Section 10], and that of the Statesman, presently under construction) has a (palindromically) periodic form.

This is the place to be reminded of that **all descriptions** of the method of (Division and) Collection in the Platonic dialogues exhbit **the cyclical, periodic nature** of the method.

In particular, the Platonic Division and Collection is not like the Aristotelian division into a genus and differentia, clearly a strictly finite division, but, having a periodic form, is an infinite division. The passage 276a1-7 is the following:

[ΝΕ. ΣΩ.]... ἀλλ' **ἡ μετὰ τοῦτο διαίρεσις**
αὖ **τίνα τρόπον** ἐγίγνετ' ἄν;
ΞΕ. **Κατὰ ταὐτὰ**

**καθ' ἅπερ ἔμπροσθεν**
**διῃρούμεθα**
τὴν ἀγελαιοτροφικὴν
πεζοῖς τε
καὶ ἀπτῆσι,
καὶ ἀμείκτοις τε
καὶ ἀκεράτοις,

[29] Continuing the previous footnote the reason given by Plato for the change from the Zeus era to the Kronos era is the need to avoid the bad, non periodic infinite, something achieved by the Logos criterion. Again Plato’s philosophic account fits well with the mathematical proof.

**τοῖς αὐτοῖς ἄν** που

**τούτοις**

**διαιρούμενοι**

**καὶ τὴν ἀγελαιοκομικὴν**

**τήν τε νῦν**

**καὶ τὴν ἐπὶ Κρόνου βασιλείαν**

**περιειληφότες** ἂν ἦμεν **ὁμοίως**

**ἐν τῷ λόγῳ.**

Translated as follows:

| | | |
|---|---|---|
| | | Younger Socrates:<br>…but how would<br>the next division be made?<br>Stranger:<br>**In the same kinds (‘kata tauta’)** |
| **as (‘kathaper’)**<br>**we divided**<br>the art of feeding herds<br>before<br>in | | |
| | those<br>that go on foot [L7] | |
| and the winged,<br>[K7] | | |
| | and<br>the unmixed breeds<br>[L9] | |
| | and the hornless,<br>[L8] | |
| | | we must **divide**<br>the art of tending herds<br>**in arts similar to**<br>**(‘tois autois’)** |
| **those (‘toutois’),** | | |
| **embracing (‘perieilephotes’) in the logos (‘en toi logoi’)**<br>both | | |
| | | the **present [Zeus]** era |
| and **the Kronos** era. | | |

| | |
|---|---|
| | The second half of the process<br>(**agelaiokomike**,<br>where the commanded herd is the human herd, held fixed, and<br>*the nature of the commanding herdsman*<br>is the variable)<br>corresponds to<br>the **Zeus** era, |
| while the first half<br>(**agelaiotrophike**, | |

<table>
<tr><td>where the nature of the commanded herd<br>is the variable)<br>corresponds to<br>the Kronos era,</td><td></td></tr>
<tr><td colspan="2">in the aneilixis myth (11.3).<br>Thus (276a) is precisely the point in the Politicus, where the method of division and collection, namely the method of philosophic periodic anthyphairesis,<br>is explicitly connected to palindromicity, since</td></tr>
<tr><td>the Kronos</td><td></td></tr>
<tr><td></td><td>and the Zeus</td></tr>
<tr><td colspan="2">eras are in palindromic relation in the myth of world ‘aneilixis’.</td></tr>
<tr><td>In this passage<br>the first half of the process,<br>with first ten division steps<br>G divided into K1, L1,<br>L1 divided into K2, L2,<br>…,<br>L9 divided into K10, L10,<br>is associated explicitly with<br>the Kronos era,</td><td></td></tr>
<tr><td></td><td>and, also explicitly, it is stated that<br>the definition of the statesman<br>must be completed,<br>as in the myth of aneilixis,<br>with its second half,<br>with last division steps, say<br>W11 divided into Z10, W10,<br>W10 divided into Z9, W9,<br>…,<br>W2 divided into Z1, W1,<br>corresponding to<br>the Zeus era,<br>palindromic to</td></tr>
<tr><td>the Kronos era,</td><td></td></tr>
<tr><td colspan="2">the whole comprising a full period, cycle.</td></tr>
</table>

(There are also the many steps L11, L12,...,W12, W11 which are not examined in detail, but are given in groups, in the form of **‘contingent causes’** (287b-289a). We will not be concerned with this part of the definition of the statesman, which probably has too some palindromic characteristics (indicated in the scheme below), as it is of no importance for our reconstruction).

In accordance with the Divisions and Collections of the *Sophistes*, this cycle must be realised by the Logos criterion.This takes the following form: The last genus of the Zeus era, denoted W1, will be divided into the art of the strategy K’1 and of the statesmanship L’1, hence the logos criterion holds K’1/L’1=K1/L1

Thus the new element in the *Politicus* Division and Collection, introduced by 11.3 and 11.4, in relation to the *Sophistes,* is not only periodicity, as in the Sophistes, but in addition the presence of palindromicity within the period; in other words we expect that the anthyphairesis will have a period, but one presenting a certain symmetry with respect to the middle of the period. In the language suggested by the myth 11.3, the division steps within the period will be divided into two parts, the first half to be the Kronos era, and the second half the Zeus era (cf. 276a discussed in (iii) and (iv) in 11e, below).

We then expect that the definition of the Statesman must exhibit the palindromicity of these two eras.We will now show that this is precisely the plan suggested and realised by Plato for defining the Statesman.

**11.5. The last ten division steps (267c-268d and 274e-277c, 289e-305e)**

The last ten steps Z10,…, Z1 in the Zeus era and the Logos Criterion in the definition of the Statesman, are described in 289e-305e.

**(i) The game (277e).**

Plato has no intention to give the Zeus steps explicitly and in succession, so these steps are not given in their correct order, in which they appear in the definition of the statesman, but in the form of a **game**, reminding of a grammatical children's game, as explicitly stated in 277e. To find the right order and the logos criterion and thus complete the definition of the statesman is going to be our task.

**(ii) The rules of the game.**

However Plato sets precise rules which enable us to complete this definition.

**Rule 1. One to One palindromic correspondence between Kj and Zj** (276a, 277e-279a).

To every division step ($K_j$, $L_j$) of the Kronos era, namely of agelaiokomike, where some type of herd $K_j$, governed and led, is **rejected,** there is a unique division step ($Z_j$, $W_j$) of the Zeus era, namely of agelaiotrophike, where some 'similar' type of governing and leadership $Z_j$ of the human herd is rejected, for j=1,..,10. The 'positive' genera Wj are **not given** (but the method of DC implies that in principle they are uniquely defined, as relative complement of $Z_j$ in $W_{j+1}$).

**Rule 2. Order of government systems** (291c-303b).

There are six types of conventional government, namely
government (I) by one, or by few, or by many, and
government (II) with laws, or without laws,
in every possible combination in each of the two criteria (I) and (II).
Thus
**T** Tyranny is government by one, without laws,
**O** Oligarchy by few without laws,
**An** Anarchy by many without laws,
**D** Democracy by many with laws,
**Ar** Aristocracy, by few with laws, and
**M** Monarchy by one with laws.

(Plato does not use the term Anarchy for the government by many without laws in the *Politicus*, but he does use the term elsewhere, e.g. in the *Laws* 942a-d).

There is a quite long discussion on these types of government in the *Politicus* 291c-303b with the following conclusions:
**(2a)**There is an ordering of these six types of government as follows:
**T>O>An>D>Ar>M**
in the sense that T is the worst and M the relatively best form of government. This ordering is so, because
(i) in the absence of government by knowledge, it is better to have laws, than to be governed without laws, and
(ii) without laws the many govern better than the few and they better than the one, while in the presence of laws the situation is precisely the reverse.
**(2b)** All these conventional types of government are to be rejected, because they are not based on true knowledge. Hence each of them (possibly more than one at a time) will be rejected at some Zj.
**(2c)** All types of government with laws are to be rejected because laws, and in fact everything written, is lifeless, as opposed to true knowledge, which has life and soul. This is strongly supported by Plato's similar rejection, in the final part of the *Phaedrus* 274c5-277a5, of the written word, in general, as an image and as not truly alive.

**Rule 3. Order of rejection** (303d-e).

Of all types of government that are to be rejected, a relatively better type of government x, rejected at stage Zj of the Zeus era, will **not be rejected before** a relatively worse type of government y, rejected at stage Zi of the Zeus era. Thus if y>x, in the sense that y is a relatively worse and x is a relatively better type of government, then i is greater than or equal to j.

**Rule 4**. **The last three arts to be rejected** (303e-305e)

The arts of Judgeship, Rhetoric and Strategy are the relatively highest arts, and are thus the last three to be rejected.

**(iii) Finding the palindromic correspondences**

We now apply these Rules, set by Plato for the purpose, in order to determine the second half of the division for the statesman

**[Z4]** It is clear by Rules 1 and 2c, considering that
in K4 lifeless herds are rejected,
that
in Z4 government by law (including **Democracy D**, **Aristocracy Ar**, **Monarchy M**) must be rejected.

**[Z3, Z6, Z7]** It is fairly clear by Rule 1 that
in Z7 the heralds (whose patron god Hermes, the herald of the Gods, is flying) (290a-c),
in Z6 the merchants and shipmasters(289e-290a), and

in Z3 the priests and diviners (290c-291a), are to be rejected as governors.

**[Z10]** Some argument is needed to convince us that the rejection of the sophists as governors in 291b-c, and again in 303b-c (of the sophists' sophists, this second time to be exact), is palindromically related to the rejection of the pigs as herd in K10 (265e-266c). Here Plato is multiply playful: in the first place the distinction of the pig from human is that it is four-footed; this is playfully described, in geometric language, as 'diameter's diameter'. We observe that the queer expression **'diameter's diameter'** has some similarity of form with the equally queer expression **'sophists' sophists'**; but we would be rather at a loss to find any more significant association; unless we note that while a true dialectician is preoccupied with 'the diameter itself' (*Politeia* 510d5-511a1), a partless Being, it is the **sophists** who, according to *Meno* 85b4, call **'diameter'** the straight line that joins one edge of a square with its opposite edge, a divisible magnitude. The meaning of this distinction is possibly the following: sophists occupy themselves with magnitudes and in general with sensibles, while the dialectician occupies himself with the idea of the diameter, something that is rather the logos of the diameter to the side of a square. In any case the *Meno* suggests that there is a connection between the diameter and the sophists; hence there is also indeed a connection between 'diameter's diameter' the four footed pig and sophists sophist, hinting, of course that the worst sophists are like pigs. We are clearly ready, granted a certain playfulness, to accept that, by Rule1, step K10, where the herd of pigs is rejected, is palindromically related with step Z10, where the sophists are rejected as governors.

**[Z8, Z9]** It now follows necessarily, by the rules 2 and 3, that the three types of conventional government without laws **T>O>An** must be associated with the three remaining free positions in the Zeus era, namely **Z9, Z8, Z5,** respectively. Namely, Tyranny must be rejected at Z9 (while K9 rejects herds of animals accepting mixed breed),

Oligarchy at Z8 (while K8 rejects herds of horned animals), and

Anarchy at Z5 (while K5 rejects herds consisting of a single animal).

It remains to see however if these necessary associations do result in palindromicity, something necessary if our interpretation is to work. For the first two of these associations, (K9, Z9) and (K8, Z8), the palindromicity follows from Plato's likening twice the two worse classes of conventional politicians with Centaurs and Satyrs (291a-b, 303c-d); thus in Z9 the Centaurs are to be rejected, representing in Plato's terminology the tyrants, in perfect palindromicity with K9, and in Z8 the Satyrs are to be rejected, representing in Plato's terminology the oligarchs, in perfect palindromicity with K8.

**[Z5]** For the last of these associations, (K5, Z5), it is not immediately obvious from the *Politicus* why this should be so, namely why Anarchy should be palindromic to the rejected herd in K4, and thus be placed in the position of the rejected part in Z4. However this is well explained by the comparison of the description of Anarchy in the *Laws* 942a5-d2, as 'kata monas dran', i.e. as a government where each person acts singly, with the description of the rejected herd in K4, as 'monotrophia', i.e. a herd consisting of a single animal Thus Anarchy (=Democracy without laws) fits well with the rejected part of Z4.

**[Z2]** We are left with the arts of Judging, Rhetoric and Strategy, for which Plato sets Rule 4, that they will be rejected last (303e-305e). Of them Judging is related to 'krinein' (305b-c), and thus Z2 is associated, by Rule 1, with Judges (by palindromicity), because K2 rejects 'krinein', and Z2 is set to reject Judges as governors and Z2 to Judges ('krinein' is the crucial word, Z2 (305b?), and K2 (260a-c).

[Z1] The art of Rhetoric is the next art rejected, as it is a practical science (304c-e), and thus Z1 is associated, by Rules 1 and 4, with Rhetoric, since K1 rejects practical sciences, and Z1 is set to reject Rhetoric.

This completes the quest for the palindromic correspondences between the first ten steps of the Kronos era with the last ten steps of the Zeus era.

**(iv) Summary of the Palindromic correspondence between the first ten steps of the Kronos era with the last ten steps of the Zeus era.**

-the rejected part of the binary division **[K10]** in the Kronos era, namely the art of governing a hornless, not mixing, four-footed (likened to **diameter's diameter**) animal herd, namely pigs, is palindromically related to
the rejected part of the binary division **[Z10]** in the Zeus era, namely the art of governing by sophists (whose practitioners are likened to **sophist's sophists**).
To see this relation note that **'diameter'** as a divisible magnitude, is opposed to 'diameter itself' (cf. *Politeia* 510d5-511a1, a partless Being), and is related to the **sophists** (as in *Meno* 85b4).

-the rejected part of the binary division **[K9]** in the Kronos era, namely the art commanding **herds of mixing with other breeds**, is palindromically related to
the rejected part of the binary division **[Z9]** in the Zeus era, namely the art of governing according to tyranny without laws, **governors likened to 'centaurs'**, **acting and governing as animals mixing with other breeds**.

-the rejected part of the binary division **[K8]** in the Kronos era, namely the art of commanding a **herd of horned animals**, is palindromically related to
the rejected part of the binary division **[Z8]** in the Zeus era, namely the art of governing according to **oligarchy** without laws, **governors likened to 'satyrs'**, **acting and governing as horned animals**.

-the rejected part of the binary division **[K7]** in the Kronos era, namely the art of commanding a **herd of flying animals,** is palindromically related to
the rejected part of the binary division **[Z7]** in the Zeus era, namely the art of **heraldry**, namely of **governing and acting according to their protector Hermes, so to speak, in a flying manner**.

-the rejected part of the binary division **[K6]** in the Kronos era, namely the art of commanding an **aquatic herd,** is palindromically related to
the rejected part of the binary division **[Z6]** in the Zeus era, namely the art of being a **merchant,** namely of **governing and acting, so to speak, in a aquatic manner**.

-the rejected part of the binary division **[K5]** in the Kronos era, namely the art of commanding **a single animal herd (‘monotrophia’)**, is palindromically related to
the rejected part of the binary division **[Z5]** in the Zeus era, namely the art of governing according to democracy without laws, namely with anarchy, described in the *Laws* 942a-d as **governing as a single person (‘kata monas dran’)**.

-the rejected part of the binary division **[K4]** in the Kronos era, namely the art of giving one’s own commands over a **lifeless herd**, is palindromically related to
the rejected part of the binary division **[Z4]** in the Zeus era, namely the art of governing **with written laws** (namely acting himself **in a lifeless manner,** according to *Phaidrus* 274c5-277a5).

-the rejected part of the binary division **[K3]** in the Kronos era, namely the art of **giving someone else’s commands,** is palindromically related to
the rejected part of the binary division **[Z3]** in the Zeus era, namely the art of priesthood (dedicated to **giving someone else’s commands)**.

-the rejected part of the binary division **[K2]** in the Kronos era, namely the art of determining, appraising (**krinein),** is palindromically related to
the rejected part of the binary division **[Z2]** in the Zeus era, namely the art of judging (which is an art of determining, appraising (**krinein)**, as is in fact explained in the *Politicus* 305b).

-the rejected part of the binary division **[K1]** in the Kronos era, namely the practical arts**,** is palindromically related to
the rejected part of the binary division **[Z1]** in the Zeus era, namely the practical art of the Rhetoric (or more precicely the part of Rhetoric dissociated from sophistry, cf. 304a) (304c-e).

**(v) Logos in the Division and Collection of the art of statesmanship**

The very last division is the division of W1 into the practical art of the strategy K’1 and the gnostical art of the statesman L’1 (304e-305a).
The **Logos criterion** consists in **equating**
**the initial logos:** practical arts **K1**/ scientific arts L1
with
**the final logos:** practical art strategy K’1/scientific art of statesmanship L’1.

**11.6. The definition of the Statesman as a philosophic palindromically periodic anthyphairesis (*Politicus* 258c-267c, 289e-305e)**

We are now ready to collect all the pieces of the Division and Collection of the Statesman.The following figure incorporates most of the information needed for the definition of the Statesman as a philosophic palindromically periodic anthyphairesis in the *Politicus.*

| Kronos era,<br>‘agelaiotrophike’ | Zeus era,<br>‘agelaiokomike’ |
|---|---|
| **[G] all arts** | |

| | | | |
|---|---|---|---|
| **(258b-259d)** | | | |
| **[K1] practical** | **[L1] scientific (gnostike)** | [W1] | [Z1] rhetoric (304c-e), |
| | | [W2] | |
| **[W1]** | | | |
| **[K’1] strategy** ***<practical>*** *(304e-305a)* | **[L’1] statesman** ***<scientific (gignoskousa)>*** | | |
| **Logos Criterion K’1/L’1=K1/L1** ***periodicity and complete knowledge of the Statesman*** | | | |
| [L1] scientific (259d-260c) | | | |
| [K2] appraise, determine **‘krinein’** | [L2] command | [W2] | [Z2] judgeship *(leaders acting on an appraising ‘krinein’ manner)* (305b-c) |
| | | [W3] (304e-305c) | |
| [L2] command (260c-261b) | | | |
| **[K3] someone else’s commands** | [L3] one’s own commands | [W3] | [Z3] priests *leaders acting on someone else’s commands)* |
| | | [W4] (290c-291a) | |
| [L3] one’s own commands (261b-d) | | | |
| **[K4] lifeless herd** | [L4] on living beings (=animals) | [W4] | [Z4] with written laws, lifeless *leaders acting in* ***lifeless*** *manner* |
| | | [W5] (291d-303b) | |
| [L4] (261d-e) on animals | | | |
| [K5] a **single** animal (‘monotrophia’) | [L5] herd of many animals | [W5] | [Z5] Democracy without laws, anarchy (‘kata monas dran’) *Leaders acting as* ***single*** |
| | | [ZW6] (291d-303b) | |
| [L5] (264b-d)<br>on many animals | | | |
| **[K6] herd of aquatic animals** | [L6] which are land animals | [W6] | [Z6] merchants (on ships)<br>*Leaders acting in **aquatic** manner* |
| | | [W7] (289e-290a) | |
| [L6] (264e-265b)<br>On many land animals | | | |
| **[K7] herd of flying animals** | [L7] which are walking | [W7] | [Z7] heralds (flying Hermes)<br>*Leaders acting in **flying** manner* |
| | | [W8] (290a-c) | |
| [L7] (265b-d)<br>On many walking land animals | | | |
| **[K8] herd of horned animals** | [L8] which are hornless | [W8] | [Z8] satyrs- (oligarchy- few, without laws)<br>*leaders acting as **horned** animals* |
| | | [W9] (291b, 303b-d) | |
| [L8] (265d-e)<br>On many hornless walking land animals | | | |
| **[K9] herd of mixing animals** | [L9] which are non-mixing | [W9] | [Z9] centaurs- tyrants (one, without laws)<br>*leaders acting as **mixing** animals* |
| | | [W10] (291a, 303b-d) | |
| [L9] (265e-266d)<br>On many non-mixing hornless walking land animals | | | |
| [K10] diameter's diameter **herd of** four-footed **pigs** | [L10] which are two-footed (human) diameter | [W10] | [Z10] sophist's sophist<br>*leaders acting as **pigs*** |
| | | [W11] (291b-c, 303b-d) | |
| | **[L10>L11>…>W11]** (287b-289a)<br>**contingent causes** | | |
| | primary | nourishment | |
| | instrument | play | |
| | receptacle | defence | |
| | vehicle | | |

The figure is to be read as follows:

**start** the division with the left column,
read from up G downwards to L10,
consisting of the ten first binary division steps of the **Kronos era**
the initial genus G is divided into K1 and L1,
L1 is divided into K2 and L2,
...,
L9 is divided into K10 and L10;

**continue** with the **contingent causes**
L10 to L11 to …to W11,
Collective steps in the division steps;

**complete** the division with the right column,
read from below W11 upwards to W1,
consisting of the ten last binary division steps of the **Zeus era**
W11 is divided into Z10 and W10,
W10 is divided into Z9 and W9,
...,
W2 is divided into Z1 and W1, and
W1 is divided into K'1 and L'1;

**end** with the **Logos criterion** K1/L1=K'1/L'1,
resulting in **periodicity** and **complete knowledge** of the Statesman.

## 12. The double Measurement in the *Politicus* 283a-287b is

## --the philosophic imitation of the Theaetetean mixture
## --the cause of Division and Collection, and
## --the Philebean Idea as mixture of infinite and finite

### 12.1. The double measurement in the *Politicus* 283a-287b is the philosophic analogue, the imitation of 'commensurability in power only'

We consider opposite entities, of the type numbers, magnitudes (lines, areas, volumes), qualities (velocities) x, y.

**There is a double 'measurement'**:
(a) x against y, and
(b) x against 'the mean(x, y)' (Greek term: 'mesos') (284e).

**Condition** on measurement **(a)**: not all opposite pairs x, y should be measured in this way; in fact if x measured against y produces a 'pleasure' (Greek term: 'hedone'), then this pair must be excluded (286d).

**Condition** on measurement**(b)**: if x, y is a pair, such that if x is measured against y then pleasure is not the product, then the second measurement, x against the mean(x, y), is a 'generation' (Greek term: 'genesis') (283d, 284b-c, d-e).

We will now decipher these terms and conditions:

**[measurement of opposites]** Measurement is about length and shortness, or excess and deficiency, or greatness and smallness, or the greater and the less, or the great and the small (283c-e), collectively of a magnitude or number or quality against its opposite (284e).

**[mean]** We interpret 'the mean(x, y)' as the philosophic analogue of the geometric mean of two numbers or magnitudes. It is remarkable that in the *Elements* there are 738 occurences of the term 'meson' (and its variants), and these exclusively and without exception refer to the 'geometric mean' (of two numbers or of two magnitudes), while NONE of them refers to the arithmetic or any other mean, or any other sense; of them, the great majority 624, occur in the Theaetetean Book X.[30]. The alternative term used 'metrios' (283e, 284a, c, d, e) has normally a meaning equivalent to 'mesos'.

**[pleasure]** That 'the dyad x, y produces **pleasure**' is equivalent, as explicitly stated in the Platonic dialogue on pleasures *Philebus* 27e, 31a, 41d, with 'the dyad x, y is an Infinite', which according to our interpretation described in Section 3a above, is equivalent with ' the dyad x, y produces an infinite'. Thus the condition on the measurement (a) is that the pair x, y does not produce an infinite, hence a finite. Thus the measurement of x against y is a (Philebean) finite.

---

[30] I owe this remark to Aliki Bassiakou [Ba].

**[generation]** In Plato there is ‘generation only’, as in the *Phaedo* 70b-72e, the *Theaetetus* 156a-157b, the *Sophistes* 247d-248a, or the *Timaeus* 27d5-29d3, 57d7-58c4, and ‘generation in substance’, as in the *Philebus* 27c, or the *Parmenides* 142c-143a or the *Politicus* 283d, 284c, d (cf. interesting comment by Proclus in *Platonic Theology* 3,91,9-24: the frequent use of words related to generation, such as ‘gignesthai’ (142d5), ‘gignetai’ (142e4), ‘genetai’ (142e6), ‘gignomenon’ (143a1), refers to ‘the progress of the intelligible multitude’ (‘peri tes proodou tou noetou plethous’). Thus ‘genesis’ and ‘generation’ does not refer only to the sensibles, but also to the intelligibles.[31] That ‘the dyad x, z is a **generation**’ implies that ‘the dyad x, z generates an infinite multitude of parts’, either as an indefinite dyad only or within a Being, namely that ‘the dyad x, z is either an Infinite or a mixture of Infinite and Finite’, hence, according to our interpretation given in Section 3, that ‘the dyad x, z produces an infinite anthyphairesis’. Thus
the condition on the measurement (b) is that the pair x, mean(x, y) produces a (Philebean) infinite.

**[The double measurement in the *Politicus* 283a-287b is the philosophic analogue of a dyad of lines incommensurable in lengh only]**
Finally in order to state the two measurement conditions (a) and (b) in the more familiar terms employed by Book X. of the *Elements*,
we choose an assumed ‘unit’ line r (as with the ‘protetheisa line’ in Book X, definition 3, and the one foot line in the Theaetetus 147-8), and,
employing Proposition II.14 of the *Elements*, we construct lines a, b, such that
$x\ r=a^2$, $y\ r=b^2$ ,
so that
mean(x, y)=a b;
The ratio x/y=the ratio $xr/yr=a^2/b^2$
hence (using the *Topics* definition, mentioned in Section 6a above),
the ratio x/ mean(x, y)=the ratio xr/m(x,y) r= the ratio $a^2/$ a b= the ratio a/b.

We are now in a position to claim that
the double measurement described is revealed to be the philosophic analogue of the mathematical condition ‘commensurability in power only’,
since
Condition (a) is equivalent to the commensurability of $a^2$, $b^2$, and
Condition (b) is equivalent to the incommensurability of a, b,
so that in conjunction they are equivalent to
the condition of commesurability in power only,
the Theaetetean mixture.

## 12.2. The double measurement as the cause of all Division and Collection, *Politicus* 283a-287b

We stay in the same passage of the *Politicus*. Plato next makes an important statement about the fundamental significance of the condition of the double measurement we have just deciphered:

[31] It is precisely the lack of realisation, by Cherniss and others, that there is intelligible ‘genesis’ as well, that has led them to mistakenly conclude that the mixtures of Infinite and Finite in the *Philebus* are not true Beings and Ideas.

(a) All arts without exception, including the statesmanship and the art of weaving (283e-284e),
(b) the method of division and collection (284e-285c), equivalently described as 'division according to kinds' (286d-287a), or as 'name and logos', sometimes shortened to 'logos' (285c-286b), by which as we saw in the *Sophistes* Platonic Beings become known to us, and
(c) the new property of 'aneilixis' (286b) , which is explained in the *Politicus* for the first time in the trilogy on knowledge,
are all caused by a double measurement as described in (12.1).
Art is a term loosely used in the *Sophistes* as well, and is never really defined formally by Plato. But we can attempt with safety to interpret this term, on the basis of its quite often use by Plato, as follows:

**[art]** A Form is knowable to humans as an Art. This is why although the Sophist is a Platonic Being, and the Statesman an Idea and a Form, nevertheless each becomes known to us in the ontological status of an Art, an inferior entity, in that divisions are more spread out and Unity is less pronounced than in the Form itself.
The term 'aneilixis' occurs in the *Politicus* only twice, in 270d and in 286b. The second occurrence (286b) has the purpose to tell us that the phenomenon of 'aneilixis' is caused in all arts by the fundamental condition of the double measurement. The first occurrence (270d) cannot have any other purpose but to DEFINE, to explain what is meant by 'aneilixis'.

### 12.3. The mixture of the Infinite and the Finite that constitutes a Platonic Idea in the *Philebus* 16c is the philosophic analogue of the Theaetetean mixture of incommensurability and commensurability; and, according to the *Theaetetus* 147d-148d Theaetetus and younger Socrates had a proof of the Proposition: every Theaetetean mixture of incommensurability and commensurability has a palindromically periodic anthyphairesis.

The two Divisions and Collections of the *Sophistes* were interpreted, in Sections 9 and 10, as the philosophic analogues of periodic anthyphairesis.
The Division and Collection of the Statesman in the *Politicus* was interpreted, in Section 11, as the philosophic analogue of palindromically periodic anthyphairesis. The cause of every Division and Collection is declared, in the *Politicus*, to be the Double measurement (Section 12.2), and the Double measurement was interpreted (Section 12.1) as the philosophic analogue of the Theaetetean mixture, namely of the dyad of lines incommensurable in length only.
Thus a remarkable consequence of our interpretation (Sections 9, 10, 11 and 12) is that the total content of the Platonic dialogue *Politicus* is the following statement:

'Every art, namely every Platonic Idea in form accessible to humans, in particular the art of the statesmanship and the art of the sophistry, possesses a a dyad of elements x, y, such that the measurement of x against y forms a philosophic finite (not a pleasure), while the measurement of x against the (geometric) mean m(x, y), forms a philosophic infinite (the generation of the mean), and precisely for this reason every art is knowable by the method of Division and Collection, which is the philosophic analogue of a palindromically periodic antrhphairesis.'

Now this philosophic statement completes the imitation 'peiro mimoumenos' that Socrates was exhorting Theaetetus to attempt, in the *Theaetetus* 147d-148d passage of the Division and Collection successfully established by Theaetetus .on a Theteatean mixture, as discussed in Setion 8.
But now we cannot fail to take note that, in fact, what has been perfectly and most closely imitated by Plato in the two dialogues Sophistes and Politicus is the modern theorem, described in Section 7, according to which ‘Every Theaetetean mixture has a palindromically periodic anthyphairesis’.

The conclusions now are the following:
--Theaetetus mathematical result described in the *Theaetetus* 147d-148d passage as Division and Collection of a power (or somewhat more generally, of a Theaetetean mixture) is precisely the modern Thyeorem (answering the question posed in Section 5.[32]
--A Platonic Idea, is the philosophic analogue of a Theaetetean mixture, namely of a (slight generalisation of) a power, and the ‘episteme’ of a Platonic Idea, obtained by the method of Division and Collectrion, is precisely the philosophic analogue of palindromically periodic anthyphairesis.
--Thus the Philebean mixture of infinite and finite is the philosophic analogue of a Theaetetean mixture.
The last two statements answer trhe questions posed in Sections 2 and 3.

**12.4. Note comparing the Philebean Finite with the Pythagorean Finite.** It should be pointed out that according to the Pythagoreans, too, the two principles of Beings is the Infinite and the Finite (cf. Philolaus, *Fragments* 1, 2, 6, and Aristotle, *Metaphysics* 987a13-28). But, while the Infinite, in both Plato and the Pythagoreans, is anthyphairetic in content, the Finite for the Pythagoreans is the Gnomon. The evolved corresponding Platonic concept is Logos, as may also be seen from a comparison of (the possibly non-genuine) *Fragment* 11 by Philolaus (knowledge is achieved by Gnomon) with the passage 546b-c of the *Politeia* on geometric number (where Gnomon is in effect replaced by Logos). But Plato in the *Philebus*, not content with such replacement, insists that the principle of Finite has anthyphairetic content as well. Plato is able to modify the original Pythagorean scheme **only because** he has now in his possession the Theaetetus palindromic periodicity theorem.

### 12.5. A note on existing interpretations of Collection and Division

The interpretation of Division and Collection presented in this paper (Sections 9, 10, 11,12) is original and runs against practically all existing interpretations.

Heidegger in [He, 1992] discusses the Division and Collection of the Angler in p.181-197, and briefly of the Sophist in . 420-422.

---

[32] Thus the view, enunciated by Knorr [Kn, 1975] that the method used by Theodorus and Theaetetus for the incommensurability proofs is not anthyphairetic, and the view, enunciated by Cherniss [Ch, 1951] and Burnyeat [Bu, 1978] that the Platonic text does not reveal the method of their incommensurability proofs are both mistaken. The method of both is shown to be anthyphairetic. Anthyphairetic reconstructions for these incommensurability proofs, mathematically impeccable but with no erguments for their historical accuracy, have been proposed by Zeuthen [Z, 1910], van der Waerden [vdW, 1950], Fowler [Fow, 1986], Kahane [Ka, 1985] (among others).

Heidegger believes that a Platonic division
-- is of finite length, ending in an indivisible species,
--is essentially the same as the Aristotelian definitions by division of genus and differentia,
and that the 'logos' in a Platonic division is identified with division.

Cornford ([Co, 1932]; [Co,1935]; [Co 1939]) thought that 'division is a downward process in which indivisible species are defined by the division of a genus with specific differences. In contrasting this approach to the characterization of Forms with earlier Platonic approaches Cornford writes: 'In a word, the Socratic method approaches the Form to be defined from below, the new method descends to it from above.' ([Mor, 1973], p.168). Collection is just the inverse of Division: he believed we could obtain the initial Genus by summing all the pieces of a Division together. According to this interpretation a Platonic Being is a One in the sense that it is the sum of its parts. Cornford and Skemp [Sk, 1952] seems to have thought that this type of collection takes place only at the beginning of the divisions. Hackforth [Ha, 1945], however, recognizes that collections take place at various places ([Mor, 1973], p.167).

Cherniss ([Ch, 1936]; [Ch, 1945]; [Ch, 1951]) assessing negatively the method of Division and Collection believed that
'the assumption of Platonic ideas is incompatible with the constitution of the species from genus and differentiae';
'The *Sophist* and *Politicus*, which have come to be considered as handbooks of *diaeresis*, show that he meant it rather to be a heuristic method, an instrument to facilitate the search for a definite idea, the distinction of that idea from other ideas, and its implications and identification, and that he did not imagine it to be a description of the "construction" of the idea, its derivation, or its constituent elements.';
'*diaeresis* appears to be only an aid to reminiscence of the ideas, .... a process the stages of which are important rather as a safeguard to insure the right direction of the search"" than as representative of necessary ingredients of the idea' ([Ch, 1945], p. 53-55).

Ryle ([Ry, 1965]; [Ry, 1966]) has an extremely low opinion of the method of Division and Collection, of the Divisions in the *Sophistes*, and of the whole *Politicus*, a 'weary dialogue'..He writes (in [Ry, 1966]):
'but until we ger to the *Sophist* we have nothing reminding us of the contribution of Linnaeus to botany; nor should we have been grateful or philosophically emlightened if we had' (*p. 136-7)*;
'the demonstrations of division in the Sophist and Politicus were intended for the tutorial benefit of beginners in the Academy',
'while the *Politicus* seems to be designed for beginners only, the Sophist is a clumsily assembled sandwich of which the bread could be of educative value only to beginners and the meat could be of value only to highly sophisticated young dialecticians'.
'So Plato may have composed the *Politicus* for the special benefit of the philosophically innocent novices who were at that moment getting their freshman's training in the ABC of thinking' (p.285).

According to a more recent assessment, by Fossheim [Fos, 2010] 'However, Ryle's judgement is all the same a difficult one to defend in face of the textual evidence. In the *Philebus* (16c), division is heralded as a divine gift, as something handed over from gods to men. Furthermore, speaking more generally, collection and division holds a prominent position in the so-called late dialogues. In short, it seems to be taken very seriously, and is unflinchingly applied to matters that are both important and difficult.'
We may add Plato's extensive expressions of admiration on the method of Division and Collection in the *Politicus* as well; and what has been said about Ryle's views on Division and Collection can also be said about the corresponding Cherniss' view.

Ackrill ([Ac, 1955]; [Ac, 1957]), in reply to Ryle, offers some sort of support of Plato's Division and Collection. But at best Plato's method for Ackrill is a primitive version of Aristotle's divisions.
He writes: 'Aristotle does indeed have, as Plato does not, a rather closely worked-out account of strict genus-species kind-ladders. But he continually speaks of divisions and of kinds and sorts ('genera' and 'species') in contexts where nobody supposes that he is, or thinks he is, doing the work of a Linnaeus.' (p.107)
'To Ackrill, its status as late Plato state-of-the-art method is corroborated by the continuity between division and the full-blown genus-species ladders developed by Aristotle.' is Fossheim's assessment [Fos, 2010].

Moravczik [Mo, 1973], p.175, in order to account for DC (*The Model of Intensional M ereology)* introduces two sorts of parts, part' and part''.
'(i) *x* is a part' of *A* =*x* has *A* as a property,+ and both *x* and *A* are Forms; furthermore, *A* does not have x as one of its properties. (E.g. an art is a part' of the Form of Art, and even more revealingly - as attested by *Sophist* 257c7ff -the parts' of the Form of Science are the various sciences.)
(ii) *x* is a part" of *A* and an 'eidos' of *A* =*A* is a Form and *x* is a kind of *A.* (E.g. acquisitive art is a kind of art, so is productive art or image-making art; these are not arts themselves like angling, sophistry, etc. but kinds of arts and thus parts" of the Form of Art.)'
'The part" relation is transitive while the part' relation is not.'
He then describes division and collection in terms of these two types of parts.
'A division consists of taking a generic Form and dividing or cutting it into a series of parts", and at the end of such a process we reach some particular part'. This part' of the generic Form is also a part' of the various parts" of the generic Form; this fact is expressed by the final definite description produced of the art in question, e.g. sophistry.'
'A collection consists of taking some parts' of *A* and showing that they are also parts' of a part" of *A.* That is to say, in a collection we show that several instances of the Form of art have some significant common feature and thus belong to a *kind* of art.'
Thus each division step, save the last one, is a division of a part'' into two parts''; in the last division step, a part'' is divided into a part'' and the final part'; and the collection is a (set-theoretic) collection of parts' into a part''.

The standard interpretations, such as given by S, Kutcharski [Ku, 1960], Miller [Mi, 1980], Scodel [Sc, 1987], Dorter [Do, 1994], Rosen [Ros, 1995], Rowe [Row, 1995], Lane [Lan, 198], Castoriadis [Ca, 1999], to name some of the most prominent, have very little in common with our mathematical analysis, in particular there is no sign or suspicion, in their interpretations, that periodic anthyphairesis plays any role.
On the contrary all existing interpretations conceive of Division and as a Linnaeus type division, strictly of finite length, essentially the same as Aristotle's divisions into genus and diferentia, and present no satisfactory explanation of Collection and Logos.

Most scholars of Plato, perhaps because of their difficulty to understand the method of Division and Collection, have argued that the theories of Knowledge that Plato describes in earlier dialogues (*Menon*, *Phaedo*, *Politeia*), namely the hypothetical method and the method by recollection, are different from the method of Diivision and Collection. However this is not correct. The method of recollection co-exists with Division and Collection in the *Phaedrus* 249b-c; the hypotheses are being divided in the *Phaedo* 107b; and the method in the *Politeia* is 'anairousa' (533c-d), quite possibly 'dividing', the hypotheses and makes the practitioner of the method, the dialectician 'sunoptikos' (537c7), namely applying Collection. The method of Recollection is naturally and simply explained in terms of the Logos Criterion for periodicity of the Collection: the Logos Criterion can be thought off as a recollection of the earlier Logos. The anthyphairetic interpretation opens the way for a satisfying unified interpretation of all the Platonic dialogues.
Ryle's low regard of Division and Collection contrasts with his high regard with the 'meat' in the *Sophistes* 'clumsy sandwich' and with the second hypothesis of the *Parmenides*. However, the whole of *Sophistes* is based on Division and Collection; and the One in the second *Parmenides* hypothesis has both 'name' and 'logos' (155d-e), another way of describing Division and Collection, and is in fact wholly based on Division and Collection and periodic anthyphairesis.

On the other side, a few students of Plato have detected in some parts of Plato's writings some direct or indirect connection with the arithmetic/geometric concept of anthyphairesis. To my knowledge, they are the following:

--Toeplitz [Toe, 1925]
stated that incommensurability must have played an important role in Plato's philosophy, but as far as I know he did not elaborate on his insight;
--Taylor [Tay, 1926] and D'Arcy W. Thompson Tho, 1929],
saw a connection between the (anthyphairetic) side and diameter numbers and the Excess and Defect described in the *Epinomis;* no wider implications were, however, realized;
--Mugler [Mug, 1948],
sensed the connection between geometric anthyphairesis via the *Theaetetus* 147-8 mathematical passage, and the *Sophistes*' Division and Collection, but in a defective way; Cherniss [Ch, 1951], in his book review 'demolished' Mugler's approach, but for the wrong reasons;
--Vuillemin [Vu, 2001],
sensed correctly that Platonic Division and Collection is related to periodic anthyphairesis, though the connection he envisaged was not the correct one. In addition; and,
--Fowler [Fow, 1986],
suggested that anthyphairesis was important in Plato's Academy, but has not explained in which way.

But these views remained decidedly marginal; researchers like Mugler and Vuillemin were unable to convince Platonists of the importance of anthyphairesis in Plato's work, while, conversely, Platonists were—and still are—unable to grasp that the Platonic method of Division and Connection describes Platonic Being. In my opinion, this double failure is due to the failure to understand the manner in which Collection turns the infinitely many parts of the anthyphairetic Division into an entity that deserves to be called One, and thus be a convincing Platonic Idea (cf. Section 14, below).

A preliminary outline of some of the basic ideas in my interpretation was given in [N, 2000]. Penrose has included some information on my interpretation in [Pe, 2004], pp. 54-59, 69.

## 13. A reconstruction of the proof of the palindromic periodicity of the anthyphairesis of commensurable in power only lines employing methods only from the Theaetetean Book X. of Euclid's *Elements*

The whole content of the *Politicus* is revealed to be nothing else but an exact philosophical imitation of the highly nontrivial theorem, stating, as mentioned in Section 7, that
'the anthyphairesis of two lines incommensurable in power only is palindromically periodic',
equivalently, in modern terms,
'the continued fraction expansion of a quadratic irrational is palindromically periodic',
a theorem thought to have been established in the modern era with fundamental contributions by some of the greatest mathematicians of the 17th and 18th century, including Euler (who gathered important experimental evidence on the periodicity of the continued fractions of particular quadratic irrationals), Lagrange (who proved the eventual periodicity), Legendre, and Galois (who characterised the purely periodic continued fractions, a result that easily implies the palindromic periodicity result). Modern proofs of this result can be found in [HW, 1938], [W, 1984] [Ka, 1985], [Fow, 1986].

The unmistakable and necessary consequence of our analysis is that, despite the lack of any explicit ancient mathematical evidence, this theorem was in possession of the Academy. Our next question is necessarily mathematical: Could the ancient Greek mathematicians, and specifically Theaetetus, obtain a proof of

**The Palindromic periodicity Theorem.** Let $\alpha>\beta$ be two lines, N a non-square number with $\alpha^2=N\ \beta^2$. Then there exist numbers $\mu_1$, $I_1$, $I_2,\ldots,I_k$, such that τηε anthyphairesis of $\alpha$, $\beta$ is palindromically periodic, and in fact equal to
$[\mu_1, \text{period}(I_1,I_2,\ldots,I_{k-1},(I_k), I_{k-1},\ldots,I_2,I_1, 2\mu_1)]$.
More generally, the same conclusion is valid for the case where $\alpha$, $\beta$ are lines commensurable in power only.

Note that both Cases actually occur. Indeed
the anthyphairesis of $\alpha,\beta$, with $\alpha^2=46\ \beta^2$, is [6, period(1,3,1,1,2,6,2,1,1,3,1,12)], and
the anthyphairesis of $\alpha,\beta$, with $\alpha^2=13\ \beta^2$, is [3, period(1,1, 1,1,6]).

Prompted by the clear, albeit in philosophical cover, statement of this theorem in the *Politicus*, as we saw in Sectionς 11-12, and even though there is no preserved ancient mathematical proof or even statement of this theorem, we next realise, that nevertheless all the tools for proving it can be found in the Theaetetean Book X. of the *Elements*.

We emphasise that our reconstruction is based solely on tools from Book X. of the Elements, but assumes a historical relevance only in connection with our discovery of this theorem hidden in philosophic language in the *Politicus*.

**13.1. The mathematical tools from the Theaetetean Book X. of the *Elements* (and the Theaetetean theory of proportion reported in the *Topics*), used for the** **reconstruction of the palindromic periodicity theorem**

We now make a complete list of the tools that will be used for the reconstruction of the proof of the palindromic periodicity theorem. These tools, with the exception of the *Topics* definition of proportion in terms of equality of anthyphairesis (which must be credited in part to Theodorus, who used its Corollary, the Logos Criterion in the reconstruction in Section V, and in greater part to Theaetetus), are all contained in the Theaetetean Book X.

The definitions of commensurability and incommensurability appear in definition 1, while the fundamental definition of commensurability in power only appears in definition 2 of Book X. of the *Elements*. Of course commensurability in power only appears explicitly in the beginning of the Platonic trilogy, in the *Theaetetus* passage 147d-148d, considered in Sections 4,5, 8), and credited to Theaetetus (and younger Socrates), and implicitly in philosophic form (as we saw in Section 12) in the *Politicus*.

The ancient concept equivalent to the modern concept of continued fraction expansion of (positive) reals is anthyphairesis of two homogeneous magnitudes; this concept is implicitly defined in the *Elements*, in the course of Propositions X.2 and X.3 (in the same way that the Euclidean algorithm for natural numbers is implicitly defined in the course of Propositions VII.1,2.

In order to fix the notation we will consider the anthyphairesis of two lines $\alpha>\beta$:

$\alpha = \mu_1\beta + \gamma_1$, with $\beta>\gamma_1$,

$\beta = I_1\gamma_1 + \gamma_2$, with $\gamma_1>\gamma_2$,

$\gamma_1 = I_2\gamma_2 + \gamma_3$, with $\gamma_2>\gamma_3$,

…

$\gamma_n = I_{n+1}\gamma_{n+1} + \gamma_{n+2}$, with $\gamma_{n+1}>\gamma_{n+2}$,

…

**Definition (Increment factors of the anthyphairetic remainders).**

Let $\alpha>\beta$ be two lines, with

$\alpha>\beta>\gamma_1>\gamma_2>\gamma_3>\ldots>\gamma_n>\gamma_{n+1}>\ldots$

the sequence of successive remainders. We define the sequence $(\varphi_n)$ of **increment factors** by the relations:

$\boldsymbol{\varphi_1=\gamma_1}$, and $\boldsymbol{\gamma_{n-1}.\varphi_n = \beta.\gamma_n}$ **for n> 1.**

Thus $\varphi_n$ is the fourth proportional of the three lines $\gamma_n$, $\beta$, and $\gamma_{n-1}$.

**Proposition 1 (Properties of the increment factors).**

Let $\alpha>\beta$ be two lines, with

$\alpha>\beta>\gamma_1>\gamma_2>\gamma_3>\ldots>\gamma_n>\gamma_{n+1}>\ldots$

the sequence of successive remainders,

$\mu_1, I_1, I_2, I_3,\ldots, I_n, I_{n+1},\ldots.$

the sequence of successive quotients, and

$\varphi_1, \varphi_2,\ldots, \varphi_n,\ldots$

the sequence of linearized remainders

of the anthyphairesis of $\alpha$, $\beta$. Then

(a) $\varphi_n<\beta$ for every n,

(b) $\varphi_n (I_n \beta+\varphi_{n+1})=\beta^2$ for every n,
(c) the anthyphairesis of β, $\varphi_n$ is equal to Anth(b, $\varphi_n$) = $[I_n, I_{n+1},\ldots]$ for every n, and
(d) [Incremental Logos Criterion] if there exist n, m with n<m, such that $\varphi_n=\varphi_m$, then the anthyphairesis of α, β is eventually periodic, and in fact equal to
$[\mu_1, I_1,\ldots, I_{n-1}$, period$(I_n, I_{n+1},\ldots, I_{m-1})]$.

**Proof**. (a), (b) routine translation of the definition of anthyphairesis of magnitudes, implicit in Proposition X.2 of the *Elements*;
(c) by the *Topika* Definition of proportion;
(d) routine translation of the Logos Criterion (established in Section 8).

**Definitions**. Let β be an assumed line, ζ, η lines with ζ>η, incommensurable, and such that $\zeta^2$, $\eta^2$ are commensurable to $\beta^2$. Then

(a) (given in Proposition X.73 of the *Elements*) the line ζ-η is called an **apotome** ('apotome') line;

(b) (given in Proposition X.36 of the *Elements*) the line ζ+η is called **line of two names** ('ek duo onomaton'); and,
(c) (implicitly given in **Propositions X.112-114** of the *Elements*) the apotome line ζ-η and the line of two names ζ+η are called **conjugate** to each other.
The lines ζ, η are called the **names** of the lines ζ-η and ζ+η.
We write $\zeta+\eta=(\zeta-\eta)^*$ and $\zeta-\eta=(\zeta+\eta)^*$.

The following, immediately verifiable, statement is nevertheless the basic link between anthyphairesis and the theory of irrational lines of Book X., introduced by Theaetetus.

**Proposition 2.** The line $\varphi_1=\alpha-\mu_1\rho$ is an apotome.

The conjugacy between the lines of two names and the apotomae is described in the fundamental Propositions X.112-114 of the *Elements*.

**Proposition X.112**. If β is a rational line, ζ is a line of two names, with names ζ1>ζ2, and η is a line such that $\zeta.\eta = \beta^2$, then η is an apotome with names $\eta_1>\eta_2$, and $\zeta_1/\eta_1=\zeta_2/\eta_2$ is a commensurable ratio.

**Proposition X.113**. If β is a rational line, ζ is an apotome, with names ζ1>ζ2, and η is a line such that $\zeta.\eta=\beta^2$, then η is a line of two names with names $\eta_1>\eta_2$, and $\zeta_1/\eta_1=\zeta_2/\eta_2$ is a commensurable ratio.

**Proposition X.114**. If ζ is a line of two names, with names $\zeta_1>\zeta_2$, η is an apotome with names $\eta_1>\eta_2$, $\zeta_1/\eta_1=\zeta_2/\eta_2$ is a commensurable ratio, and $\eta.\zeta=\beta^2$, then β is a rational line.

In these three Propositions Theaetetus develops the conjugacy of the two types of lines (apotome, of two names) he has introduced; this conjugacy plays in effect the same role that conjugation plays in complex numbers: it is useful in computing inverses, a computation that is vital for the continued fraction expansions. The reconstruction of the proof of the palindromic periodicity theorem uses solely the concepts of the lines apotome and of two names, and their conjugacy as stated in Propositions 112-14 of Book X. of the *Elements*.

**13.2. The tools of Book X are sufficient to carry Euler-type computations on anthyphairetic computations on ‘powers’.**

Before embarking on the proof of the Periodicity and Palindromicity Theorems it is instructive to get a feeling of the importance of the conjugacy of the binomial and the apotome lines in the anthyphairesis of quadratic irrationals, by translating the calculations performed by Euler, in computing the continued fraction expansion of his favorite quadratic irrational, into the language of the tools of Book X. of the *Elements*, i.e. considering them from Theaetetus viewpoint (cf. [Fow, 1986], pp. 322-323).
Let us note that the importance of these propositions has not really been realized, e.g. by Heiberg, as it appears from Heath’s comments in [Hea, 1926], p. 246, which however are in the right direction:
**“**Heiberg considers that this proposition and the succeeding ones are interpolated, though the interpolation must have taken place before Theon’s time. His argument is that X.112-115 are nowhere used.

It seems to me that [Propositions X.112-114] that they *are* so connected…. But X.112-114 show us how either of them (binomial straight line, apotome) can be used to rationalize the other, thus giving what is surely an important relation between them.’

Weil [W, 1984] (p.15), evidently not aware of Plato’s association with periodic anthyphairesis, writes: ‘Not only is Euclid himself well aware of the relation
$(\sqrt{r}+\sqrt{s})(\sqrt{r}-\sqrt{s})=r-s$
but even the identity
$1/(\sqrt{s}+\sqrt{r}) = \sqrt{r}/(r-s) - \sqrt{s}/(r-s)$
may be regarded as the essential content of prop. 112 of that book. Unfortunately, Euclid's motivation in Book X seems to have been the wish to construct a general framework for the theory of regular polygons and polyhedra, and not, as modern mathematicians would have it, an algebraic theory of quadratic fields. So we are left to speculate idly whether, in antiquity or later, identities involving square roots may not have been used, at least heuristically, in arithmetical work.’

| **Euler [E, 1765] calculations on continued fraction expansions of quadratic irrationals** | **Corresponding calculations on anthyphairetic expansions of ‘powers’, employing the tools of Book X. of Euclid’s *Elements*** |
|---|---|
| $\sqrt{54}$ | $\alpha^2 = 54\beta^2$ |
| First division step | |
| $\mu_1$=integral part of $\sqrt{54}=7$ | $\mu_1=7$ |
| $x_1=\sqrt{54}-7$ | [first anthyphairetic remainder $\gamma_1=\varphi_1$] $\varphi_1=\alpha-7\beta$. The first anthyphairetic remainder $\varphi_1$ is an **apotome** |
| Second division step | |
| $y_1= 1/(\sqrt{54}-7)=(\sqrt{54}+7)/5$ | [inversion and conjugation] |
| | Let $\psi_1$ be the **line inverse** to the apotome $\varphi_1$, namely the line that satisfies<br>$\varphi_1.\psi_1= \beta^2$.<br>The fact that $\varphi_1$ is an apotome allows for the exact computation of the inverse line $\psi_1$.<br>Let $\varphi_1$* be the **line of two names conjugate** to $\varphi_1$, namely $\varphi_1^*=\alpha -+7\beta$.<br>By conjugation (**Proposition X.112**)<br>$\varphi_1.\varphi_1^*=\alpha^2-49\beta^2=54b^2-49\beta^2=5\beta^2$.<br>Comparison of $\psi_1$ and $\varphi_1$* yields immediately<br>$5\ \psi_1 =\alpha -+7\beta$.<br>Thus $\psi_1$ is a line of two names. |
| $I_1$=integral part of 14/5=2 | [inversion and anthyphairesis]<br>We compare<br>$\varphi_1.\psi_1=\beta^2$ (inversion), and<br>$\varphi_1.(I_1\ \beta+\varphi_2)= \beta^2$ (increment factor).<br>hence $\psi_1=I_1\ \beta+\varphi_2$.<br>Since $\varphi_2< \beta$, and $5\ \psi_1 =\alpha -+7\beta$,<br>It follows that<br>**$I_1$** is the **integral part** of 14/5,<br>namely $I_1$=2. |
| $x_2= (\sqrt{54}-3)/5$ | [increment factor $\varphi_2$]<br>hence $\alpha+ \mu\beta=5(2\ \beta+\varphi_2)$,<br>hence $5\varphi_2=\alpha- 3\beta$,<br>$\varphi_2$ is an apotome |
| Third division step | |
| $y_2=5/(\sqrt{54}-3) =5(\sqrt{54}+3)/45 = (\sqrt{54}+3)/9$ | [inversion and conjugation]<br>$\varphi_2.\psi_2=\beta^2$, thus $9\psi_2=\alpha+ 3\beta$<br>thus $\psi_2$ is the line of two names inverse to the apotome $\varphi_2$,<br>found by employing **Proposition X.112** |
| $I_2$= integral part of 10/9=1 | [inversion and anthyphairesis]<br>We compare<br>$\varphi_2.\psi_2=\beta^2$ (inversion), and<br>$\varphi_2.(I_2\ \beta+\varphi_3)= \beta^2$ (increment factor).<br>hence $\psi_2=I_2\ \beta+\varphi_3$.<br>Since $\varphi_3< \beta$, and $9\psi_2=\alpha+ 3\beta$<br>**$I_2$** is the integral part of 10/9,<br>namely $I_2$=1. |
| $x_3=(\sqrt{54}- 6)/9$ | [increment factor $\varphi_3$]<br>hence $\alpha+ 3\beta=9(\beta+\varphi_3)$,<br>hence $9\varphi_3=\alpha- 6\beta$,<br>The increment factor $\varphi_3$ is an apotome |
| Fourth division step | |
| $y3=9/(\sqrt{54}-6)=9(\sqrt{54}+6)/18=(\sqrt{54}+6)/2$ | [inversion and conjugation]<br>$\varphi_3.\psi_3=\beta^2$, thus $2\psi_3=\alpha+ 6\beta$<br>thus $\psi_3$ is the line of two names inverse to the apotome $\varphi_3$, |
| | found by employing **Proposition X.112** |
| $I_3=6$ | [inversion and anthyphairesis]<br>We compare<br>$\varphi_3.\psi_3=\beta^2$ (inversion), and<br>$\varphi_3.(I_3\ \beta+\varphi_4)=\beta^2$ (increment factor).<br>hence $\psi_3=I_3\ \beta+\varphi_4$.<br>Since $\varphi_4<\beta$,<br>$\mathbf{I_3}$ is the integral part of 13/2, namely<br>$I_3=6$. |
| $x_4=(\sqrt{54}-6)/2$ | [increment factor $\varphi_4$]<br>hence $\alpha+6\beta=2(6\beta+\varphi_4)$,<br>hence $2\varphi_4=\alpha-6\beta$,<br>The increment factor $\varphi_4$ is an apotome |
| Fifth division step | |
| $y_4=2/(\sqrt{54}-6)=2(\sqrt{54}+6)/18=(\sqrt{54}+6)/9$ | [inversion and conjugation]<br>$\varphi_4.\psi_4=\beta^2$, thus $9\psi_4=\alpha+6\beta$<br>thus $\psi_4$ is the line of two names<br>inverse to the apotome $\varphi_4$,<br>found by employing **Proposition X.112** |
| $I_4=1$ | [inversion and anthyphairesis]<br>We compare<br>$\varphi_4.\psi_4=\beta^2$ (inversion), and<br>$\varphi_4.(I_4\ \beta+\varphi_5)=\beta^2$ (increment factor).<br>hence $\psi_4=I_4\ \beta+\varphi_5$.<br>Since $\varphi_5<\beta$,<br>$\mathbf{I_4}$ is the integral part of 13/9, namely<br>$I_4=1$.<br>Thus $\psi_4=\beta+\varphi_5$. |
| $x_5=(\sqrt{54}-3)/9$ | [increment factor $\varphi_5$]<br>hence $\alpha+6\beta=9(\beta+\varphi_5)$,<br>hence $9\varphi_5=\alpha-3\beta$,<br>The increment factor $\varphi_5$ is an apotome |
| Sixth division step | |
| $y_5=9/(\sqrt{54}-3)=9(\sqrt{54}+3)/45=(\sqrt{54}+3)/5$ | [inversion and conjugation]<br>$\varphi_5.\psi_5=\beta^2$, thus $5\psi_5=\alpha+3\beta$<br>thus $\psi_5$ is the line of two names<br>inverse to the apotome $\varphi_5$,<br>found by employing **Proposition X.112** |
| $I_5=2$ | [inversion and anthyphairesis]<br>We compare<br>$\varphi_5.\psi_5=\beta^2$ (inversion), and<br>$\varphi_5.(I_5\ \beta+\varphi_6)=\beta^2$ (increment factor).<br>hence $\psi_5=I_5\ \beta+\varphi_6$.<br>Since $\varphi_6<\beta$,<br>$\mathbf{I_5}$ is the integral part of 10/5, namely<br>$I_5=2$.<br>Thus $\psi_5=2\beta+\varphi_6$. |
| $x_6=(\sqrt{54}-7)/5$ | [increment factor $\varphi_6$] hence $\alpha+3\beta=5(2\beta+\varphi_6)$, |
| | hence $5\varphi_6=\alpha- 7\beta$,<br>The increment factor $\varphi_6$ is an apotome |
| Seventh division step | |
| $y_6=5/(\sqrt{54}-7)=5(\sqrt{54}+7)/5=\sqrt{54}+7$ | [inversion and conjugation]<br>$\varphi_6.\psi_6=\beta^2$, thus $\psi_6=\alpha+ 7\beta$<br>thus $y_6$ is the line of two names<br>inverse to the apotome $\varphi_6$,<br>found by employing **Proposition X.112** |
| $I_6=14$ | [inversion and anthyphairesis]<br>We compare<br>$\varphi_6.\psi_6=\beta^2$ (conjugacy), and<br>$\varphi_6.(I_6\beta+\varphi_7)= \beta^2$ (increment factor).<br>hence $\psi_6=I_6 \beta+\varphi_7$.<br>Since $\varphi_7< \beta$,<br>$\mathbf{I_6}$ is the integral part of 14, namely $I_6=14$.<br>Thus $\psi_6=14\beta+\varphi_7$. |
| $x_7=\sqrt{54}-7$ | [increment factor $\varphi_7$] hence $\alpha+ 7\beta=14\beta+\varphi_7$,<br>hence $\varphi_7=\alpha- 7\beta$,<br>The increment factor $\varphi_7$ is an apotome |
| Logos Criterion | |
| $x_7=x_1$, periodicity | $\boldsymbol{\varphi_7=\varphi_1}$, periodicity<br>by the incremental **Logos criterion** |
| Hence the continued fraction expansion of $\sqrt{54}$ is equal to<br>[7, **period**(2,1,6,1,2,14)]. | Hence the anthyphairesis of $\alpha$, $\beta$, with $\alpha^2= 54\beta^2$, is equal to<br>[7, **period**(2,1,6,1,2,14)]. |

**13.3. The proof of the Periodicity Theorem with the tools of Book X of the *Elements***

We now turn to the proof of the Theorem, and present first the proof of the weaker

**The periodicity theorem (PT)**.
If N is a non-square natural number, and $\alpha$ and $\beta$ straight lines, such that $\alpha^2=N\beta^2$, denoting by
$\mu_1, I_1, I_2,\dots, I_n,\dots$
the sequence of successive quotients, by
$\varphi_1, \varphi_2,\dots, \varphi_n,\dots$
the sequence of successive increment factors.
of the anthyphairesis of $\alpha$ with respect to $\beta$,
then there are two sequences
$\lambda_1(=1), \lambda_2,\dots, \lambda_k,\dots$ and $\mu_1, \mu_2,\dots, \mu_k,\dots$
of natural numbers, such that for every natural number k>1
($a_k$) $\lambda_k\lambda_{k-1}=N- \mu_{k-1}^2$,
($b_k$) $\lambda_k (I_{k-1} \beta +\varphi_k)= \lambda_{k-1} (\varphi_{k-1})^*$,
($c_k$) $\mu_k+\mu_{k-1}=I_{k-1}\lambda_k$,
hence $\lambda_k$ is a natural number smaller than N, and $\mu_k$ is a natural number whose square is smaller than N, and
($d_k$) the line $\varphi_k$ is an **apotome** and in fact $\lambda_k \varphi_k=\alpha-\mu_k \beta$; and
the anthyphairesis of $\alpha$ and $\beta$ is eventually periodic.

**Proof of the periodicity theorem.**
We proceed by mathematical induction.
We prove ($a_2$), ($b_2$),($c_2$), ($d_2$) and
assuming ($a_k$), ($b_k$),($c_k$), ($d_k$) we prove ($a_{k+1}$), ($b_{k+1}$),($c_{k+1}$), ($d_{k+1}$).

The proof of ($a_{k+1}$) consists in the elementary verification, using (($d_k$), ($a_k$), ($c_k$)), that the number $\lambda_{k+1}$ defined by $\lambda_{k-1}+I_{k-1}(\mu_{k-1}-\mu_k)$ is indeed a natural number.

Proof of ($b_{k+1}$): the equality ($d_k$) $\lambda_k\varphi_k=\alpha-\mu_k\beta$ implies that $\varphi_k$ is an apotome line. We consider the line of two names $\lambda_k(\varphi_k)^*=\alpha+\mu_k\beta$ conjugate to it; by ($a_{k+1}$) and **Proposition X.114**, we obtain that

$$\lambda_k^2\varphi_k(\varphi_k)^*=(\alpha-\mu_k\beta)(\alpha+\mu_k\beta)=\alpha^2-\mu_k^2\beta^2=\lambda_{k+1}\lambda_k\beta^2,$$

whence, dividing both sides by $\lambda_k$, that

$$\lambda_k\varphi_k(\varphi_k)^*=\lambda_{k+1}\beta^2. \qquad \text{(A)}$$

On the other hand, from the anthyphairetic relation (Proposition 1b) we obtain $\varphi_k(I_k\beta+\varphi_{k+1})=\beta^2$, hence

$$\lambda_{k+1}\varphi_k(I_k\beta+\varphi_{k+1})=\lambda_{k+1}\beta^2. \qquad \text{(B)}$$

Comparing (A) and (B), we get that $\lambda_{k+1}\varphi_k(I_k\beta+\varphi_{k+1})=\lambda_k\varphi_k(\varphi_k)^*$, hence

$$\lambda_{k+1}(I_k\beta+\varphi_{k+1})=\lambda_k(\varphi_k)^*.$$

Assuming that ($b_{k+1}$) has been proved:
the proof of ($c_{k+1}$) consists in the elementary verification, using (($b_{k+1}$) and Proposition 1a), that the number $\mu_{k+1}$ defined by $I_k\lambda_{k+1}-\mu_k$ is indeed a natural number; and,
the proof of ($d_{k+1}$) consists in an elementary verification, using ($b_{k+1}$) and ($c_{k+1}$).

The similar and simpler proof of ($b_2$), relying on the conjugacy of Book X. (**Proposition X.114**), is left for the reader.

The proof of ($a_k$), ($b_k$),($c_k$), ($d_k$) is now complete.

From ($a_k$), ($b_k$),($c_k$), ($d_k$) just proved, and from the fact that the set of pairs ($\lambda$, $\mu$) of natural numbers, with $\lambda$ bounded above by N and $\mu$ by $\mu_1$ is finite, it follows, by an application of **the pigeonhole principle**, that there are indices k and s, with k<s, such that $\varphi_k=\varphi_s$. The eventual periodicity now follows from the incremental Logos criterion (Proposition 1d).

### 13.4. The proof of the palindromic periodicity theorem with the tools of Book X of the *Elements*

We next proceed with the proof of the full theorem on palindromic periodicity, described in the *Politicus*.

**Definition (Inverse of the conjugate ($\varphi_n$)*).** Define the line $\omega_n$ by the equality $(\varphi_n)^*\omega_n=\beta^2$.
Since $(\varphi_n)^*$ is a line of two names, it follows from Book X. conjugacy (**Proposition X.112)** that $\omega_n$ is an apotome. In fact we can obtain more detailed information.

**Proposition 3 (Properties of $\omega_n$).**
(a) The line $\omega_n$ is an apotome, and in fact

$\lambda_{n+1}\,\omega_n=\alpha-\mu_n\,\beta$;
(b) $\omega_n<\beta$;
(c) $\omega_1(\varphi_1+2\mu_1\,\beta)=\beta^2$; and,
(d) $\omega_{n+1}\,(I_n\,\beta+\omega_n)=\beta^2$.

**Proof**. (a) by the definition of $\omega_n$ and PT(=Periodicity Theorem) ($a_n$), ($d_n$);

(b) by the definition of $\omega_n$ and PT ($c_n$);

(c) by the definition of $\omega_1$ and $\varphi_1$; and,
(d) by the definition of $\omega_n$ and PT ($d_{n+2}$).

The details are left for the reader.

**Proof of the palindromic periodicity theorem.**

Order the elements of the two sequences ($\varphi_n$) and ($\omega_n$) in one as follows:
$\varphi_1, \omega_1, \varphi_2, \omega_2,\ldots, \varphi_n, \omega_n,\ldots$ .
From **the pigeonhole principle**, using the Periodicity Theorem and Proposition 3a, there exists a first element in the sequence so ordered, such that
(i) all the previous elements of the sequence are pairwise distinct, and
(ii) the element in question coincides with a previous element.

There are two cases to examine:
**Case I**: The elements $\varphi_1, \omega_1, \varphi_2, \omega_2,\ldots, \varphi_{k-1}, \omega_{k-1}$ are pairwise distinct and the element $\varphi_k$ coincides with a previous element. Then we claim that $\omega_{k-1}=\varphi_k$ .
**Case II**: The elements $\varphi_1, \omega_1, \varphi_2, \omega_2,\ldots, \varphi_{k-1}, \omega_{k-1}, \varphi_k$ are pairwise distinct and the element $\omega_k$ coincides with a previous one. Then we claim that $\varphi_k=\omega_k$.

**Proof for Case I**. There are three subcases:

(a) if the apotome line $\varphi_k$ coincides with the apotome line $\varphi_1$, then we have contradiction [Proof similar to (b) below]

(b) if the apotome line $\varphi_k$ coincides with the apotome line $\varphi_j$ for some j with 1<j<k, then we have contradiction.
[we are to obtain a contradiction from (b).
From the supposition it follows that the two conjugate lines of two names are equal, i.e. that $(\varphi_k)^*=(\varphi_l)^*$, and since $(\varphi_k)^*\,\omega_k=(\varphi_j)^*\,\omega_j=\beta^2$ (by the definition of $(\varphi_k)^*$), we obtain that $\omega_k=\omega_j$. But we have $\omega_j\,(I_{j-1}\,\beta+\omega_{j-1})=\beta^2$and $\omega_k\,(I_{k-1}\,\beta+\omega_{k-1})=\beta^2$ (both by Proposition 3d). Hence $I_{k-1}\,\beta+\omega_{k-1}= I_{j-1}\,\beta+\omega_{j-1}$; since (on account of Proposition 3b) both sides of the equation represent the anthyphairetic division of the same whole line with respect to the same line $\beta$, it follows that the remainders are equal, i.e. $\omega_{k-1}=\omega_{j-1}$, a contradiction, since the elements $\varphi_1, \omega_1, \varphi_2, \omega_2,\ldots, \varphi_{k-1}, \omega_{k-1}$ have been assumed to be pairwise different. Hence $\varphi_k$ is not equal to $\varphi_j$ for some j with 1<j<k.]

(c) if the apotome line $\varphi_k$ coincides with the apotome line $\omega_j$, for some j with j<k-1, then we have contradiction.[Proof similar to (b) above].

Therefore, since $\varphi_k$ is equal with one of the elements $\varphi_1, \omega_1, \varphi_2, \omega_2,\ldots, \varphi_{k-1}, \omega_{k-1}$, the only possibility left, by (a), (b), and (c), is to have $\omega_{k-1}=\varphi_k$.

**Case II** is dealt with similarly and we obtain $\varphi_k = \omega_k$.

We now complete the proof by showing the palindromicity in **Case I**.

From the fact, already proved, that $\omega_{k-1} = \varphi_k$, and the relations $\varphi_k\ (I_k\beta+\varphi_{k+1})=\beta^2$ (Proposition 1b) and $\omega_{k-1}\ (I_{k-2}\ \beta+\omega_{k-2})=\beta^2$ (Proposition 3d), we deduce that
$I_k\ \beta+\varphi_{k+1} = I_{k-2}\ \beta+\omega_{k-2}$;
since (on account of Propositions 1a and 3b) both sides of this equation represent the anthyphairetic division of the same whole line with respect to the same line $\beta$, it follows that the two remainders and the two quotients are equal, i.e. $I_k=I_{k-2}$ and
$\varphi_{k+1} = \omega_{k-2}$ .
Continuing in the same way we obtain that $I_{k+1}=I_{k-3}$ and $\varphi_{k+2} = \omega_{k\,3}$. Still continuing in the same way we obtain that $I_1=I_{2k-3}$ and $\varphi_{2k+2} = \omega_1$. But we have that
$\varphi_{2k-2}(I_{2k-2}\ \beta+\varphi_{2k-1})=\beta^2$ (Proposition 1b) and $\omega_1(2\mu_1\beta+\varphi_1)=\beta^2$ (Proposition 3d); hence $I_{2k-2}\ \beta+\varphi_{2k-1} = 2\mu_1\beta+\varphi_1$. Since (on account of Proposition 1a) both sides of this equation represent the anthyphairetic division of the same whole line with respect to the same line $\beta$, it follows that the two remainders and the two quotients are equal, i.e. $I_{2k-2}=2\mu_1$ and $\varphi_{2k-1} = \varphi_1$.
From the equality $\varphi_{2k-1} = \varphi_1$ (the incremenal Logos criterion (Proposition 1d)), we deduce the periodicity of the anthyphairesis, and from the equalities
$I_{2k-2}=2\mu_1$, and $I_1=I_{2k-3}$, $I_2=I_{2k-4},\ldots, I_k=I_{k-2}$, and $I_{k+1}=I_{k-3}$,
we deduce the palindromicity in the period of the anthyphairesis, as required.

**Case II** is dealt with similarly.

The general commensurability in power only requires only minor modifications left to the reader. The proof of the theorem is complete.

**13.5. Plato's argument for going from the Kronos era to the Zeus era (*Politicus* 272d-e) appears to be related to the pigeonhole argument employed in the proof.**

It will be seen that in the reconstructed proof of the palindromic periodicity theorem, the pigeonhole principle is crucially employed precisely at the point of passing from the first semiperiod to the second and palindromically symmetric to the first. This of course corresponds precisely to the point of passing from the Kronos era to the Zeus era. As noted in Section 11.3, and footnote 22, Plato, in explaining why there comes a time when there is necessarily a change from the Kronos era to the Zeus era (in the *Politicus* 272d-e), gives an argument that seems reminiscent of a pigeonhole argument:
'For when the time of all these [states] was **completed** ('eteleiothe') and the change was to take place, and ***all the earth-born genus had already been used up ('aneloto'), since every soul falling into the earth had fulfilled ('apodedokuias') all the generations, as many seeds as prescribed to it,*** then the heldsman of the universe dropped the tiller and withdrew to his place of outlook, and fate and innate desire made the universe turn backwards.'[33]

---

[33] ἐπειδὴ γὰρ πάντων τούτων χρόνος ἐτελεώθη καὶ μεταβολὴν ἔδει γίγνεσθαι καὶ δὴ καὶ τὸ γήινον ἤδη πᾶν ἀνήλωτο γένος, πάσας ἑκάστης τῆς ψυχῆς τὰς γενέσεις ἀποδεδωκυίας, ὅσα ἦν ἑκάστῃ προσταχθὲν τοσαῦτα εἰς γῆν σπέρματα πεσούσης

Thus change for a soul from the Kronos earthborn era to the Zeus era becomes necessary when the (definitely finite) number of generations of the soul in this era have all been used up and the seeds prescribed to the soul in this era have been; this is precisely the pigeonhole argument.

### 13.6. The purpose of Book X of the *Elements*

Book X. is thus revealed to have a most exciting purpose: supplying the tools for the statement and proof of the palindromic periodicity of the anthyphairesis of the quadratic irrationals. This purpose was hitherto unsuspected, as can be seen by opinions expressed by authors who have studied Book X., such as:
Taisbak ([Tai, 1982], p.58), who must be credited with demystifying the mathematical structure of Book X., writes: "We are prepared to face the possibility that there was no other point than to entertain us with good logic";
Knorr ([Kn, 1983], p. 60), who writes: "The true merit of Book X, and I believe it is no small one, lies in its being a unique specimen of a fully elaborated deductive system of the sort that ancient philosophies consistently prized";
Mueller ([Mu, 1981], pp. 270-271), who writes: "One would, of course, prefer an explanation that invoked a clear mathematical goal intelligible to us in terms of our own notions of mathematics and which, under analysis, would lead univocally to the reasoning in book X. Unfortunately, book X has never been explicated successfully in this way nor does it appear amenable to explication of this sort. Rather, book X appears to be expedient for dealing with a particular problem and at the same time a mathematical blind alley";
van dr Waerden ([vdW, 1950], p. 172), who writes: " Book X does not make easy reading…The author succeeded admirably in hiding his line of thought".

On the other hand, the 'fine structure' of Book X., involving the detailed study of the six lines of two names and the six apotomai has no role in the proof of the palindromic periodicity theorem, and its purpose is still an open question. It will be analysed in a forthcoming work by Negrepontis and Brokou [NB, ?].

## 14. A Platonic Idea is an Anthyphairetic Self-Similar One

There is a final question that must be answered, in order to consider the interpretation given in the present work fully satisfying. The question is: why would Plato want to place periodic anthyphairesis at the very center of his philosophic system? What was he expecting from periodic anthyphairesis?
We already saw that by identifying Platonic Ideas with the philosophic analogue of a Theaetetean mixture he succeeded in making the Ideas knowable to humans by the method of Division and Collection.
The second advantage of periodic anthyphairesis is that it possesses an infinite sequence of anthyphairetic approximations, namely of rational approximations whose finite anthyphairesis is a finite initial segment of the infinite anthyphairesis of the Platonic Idea.Thus the anthyphairetic model of an Idea provides approximations,the sensible, that may credibly be described as copies or idols of, or participants into, the original Idea. Sensibles are expanded by Plato mostly in the *Timaeus*. As already noted we are not going to deal with sensible in the present work.

But the main property of a Platonic Idea that Plato obtains from the model of periodic anthyphairesis is its Oneness. In fact an entity possessing periodic anthyphairesis is a One, a Monad, in the best possible sense, that of self-similarity. We will then describe in some detail the precise sense in which a Platonic Idea or Being may be regarded as a One, as promised in Section 1.

### 14.1. What kind of a One for the Platonic Idea?

Plato distinguishes two states of Beings in the *Sophistes* 255c-d: the more exalted ‘Being itself’ and a more lowly Being ‘with respect to’ (‘pros ti’), a relative Being, a Being in the form of a ratio, an Idea, according to our foregoing analysis. Only the latter is accessible and knowable to the human intellect, and we are here concerned exclusively with this Being. It is clear that, for Plato, the principal property of a Platonic Being—variously described as a One or a Unity, as simple, partless or indivisible—is precisely its Oneness. But what kind of Oneness?

### 14.2. The types of One rejected for a Platonic Being

*(1) A Platonic Being is not the absolutely partless One.*

The Absolute One, the One without any parts whatsoever, the really partless and indivisible One which corresponds to the One of the first hypothesis (137c-142a) in the *Parmenides*, is explicitly rejected as Being (especially in 141e3-142a1). In this rejection, Plato follows Zeno, as in *Fragment B2*.

*(2) A Platonic Being is not the cumulative One.*

Cornford ([Co, 1932]; [Co,1935]; [Co,1939]) thought that Collection is just the inverse of Division: he believed we could obtain the initial Genus by summing all the pieces of a Division together (cf. Section 12.5). But there are serious problems with this interpretation of Collection, and the corresponding interpretation of a Platonic

Being as a One that is simply the sum of an infinity of parts, a summed into a totality. To begin with, there is no periodicity involved in such a Collection, although, as we note in Section 8, every Platonic description of Collection stresses its periodic nature. Secondly, the type of Collection suggested by Cornford is explicitly rejected both in the *Parmenides* 128e-130a (especially 129c4-d6 given in 5(c), below) and in the *Philebus* 14d-e, as describing sensible and not intelligible entities:

*'and this sort of thing also should be disregarded, when a man in his talk dividing the members, and at the same time the parts, of anything, acknowledges that they all collectively are that one thing'*[34] (Philebus 14e).

### 14.3. A Platonic Idea is described as a self-similar One in the *Sophistes* and in the second hypothesis of the *Parmenides*

At first sight, infinite divisibility seems to make a mockery of—and to run against—any decent notion of One, as it produces an infinite number of parts, namely the remainders at each stage of the division process. But as it is made clear in the *Sophistes* 244d-245b, 257c, 258e, the *Parmenides* 129c-d, and the *Parmenides* Second Hypothesis, divisibility is no obstacle for Oneness, specifically for self-similar Oneness in which every part is the same as the whole:

*'But yet nothing hinders ('ouden apokoluei') that which has parts ('to memerismenon') from possessing the attribute of unity ('to pathos …tou henos') in all its parts ('en tois meresin pasin'), and in this way every Being and Whole ('pan te on kai holon') to be One ('hen')'*[19] (*Sophistes* 245a1-3);

*'but what is surprising if someone shall show that I am one and many ('hen kai polla')? When he wishes to show that I am many, he says that my right side is one thing and my left another, that my front is different from my back, and my upper body in like manner different from my lower; for I suppose I have a share of multitude. To show that I am one, he'll say I am one man among the seven of us, since I also have a share of the one. So he shows both are true. Now, if someone should undertake to show that sticks and stones and things like that are many, and the same things one, we'll grant that he has proved that something is many and one ('polla kai hen'), but not that the one is many ('to hen polla') or the many one ('ta polla hen'): he has said nothing out of the ordinary, but a thing on which we all agree.'* (*Parmenides* 129c4-d6);

*'there are infinitely many parts of Being' (Parmenides 144b6-7), 'to each of the parts of Being befits the One' (ibid. 144c6),'neither the being is lacking in relation to the one, nor the one in relation to the being, but they are equalized being two for ever in all ways'* (*ibid.* 144e1-3).[35]

[34] Based on Plato, *Statesman, Philebus*, translated by H.N. Fowler, Loeb Classical Library, Cambridge, Mass., 1925.

[35] translation based on Allen [Al. 1997].

### 14.4. A Platonic Idea is an Anthyphairetic Self-similar One

*(1) Anthyphairetic Division into an infinite multitude of parts.*

An Idea possesses, in analogy to a mathematical ratio of lines incommensurable in length only, two opposite to each other parts called powers in the *Theaetetus* 156a6-b1—$B_1$ and $A_1$—one of which—say $A_1$—acting and the other—$B_1$—being acted on (*Sophistes* 247d-e); this is the Platonic entity, called by Aristotle an 'indefinite dyad' in, for example, the *Metaphysics* 1082a13-15, 1083b35-36.
Examples of Platonic 'indefinite dyads' defining a Being are:
{$A_1$ beautiful, $B_1$ not beautiful}(*Sophistes* 257d-e);
{$A_1$ great, $B_1$ not great} and {$A_1$ just, $B_1$ not just}(*ibid*. 258a);
{being $A_1$, non-being $B_1$} (*ibid*. 258b-c);
{$A_1$ self-restraint and $B_1$ bravery}, the indefinite dyad of the Form Virtue (mentioned and analysed in the final section of the *Politicus* 305e-311c); and
{$B_1$ One, $A_1$ Being}, the indefinite dyad of the Form One Being (in the second hypothesis of the *Parmenides* 142b1-143a3).
Indefinite dyads without an explicit reference to a Platonic Idea, such as
{cold, warm}, {fast, slow}, {more, less},
appear in the *Philebus* 23b5-25e3 and are instances of the Infinite ('apeiron'), interpreted in Section 3, as a philosophic form of infinite anthyphairesis.
The initial indefinite dyad of the Ideas the Angler and the Sophist is, in each case, the pair of the opposite species-parts {$B_1$, $A_1$} produced at the first stage by dividing the initial genus G.
As shown in the second hypothesis of the *Parmenides* 142b1-143a3[36], for the dyad {One, Being}, the initial infinite dyad {$B_1$, $A_1$} produces, by anthyphairetic division, the infinite sequence of parts-species
$A, B, B_1, A_1, B_2, A_2, \ldots, B_k, A_k, B_{k+1}, A_{k+1}, \ldots$.
In the abbreviated form of the Division, as it appears in the Divisions in the *Sophistes* and in the *Politicus*, the parts of a Platonic Being are produced by dividing the initial Genus G at the first stage into two species $B_1$ and $A_1$, one of which—the one which is a true predicate to the Platonic Being, $A_1$—turns into a genus and is divided again at the second stage into two opposite species-parts, say $B_2$ and $A_2$. Because of Logos-periodicity, the Division continues ad infinitum, producing the infinite sequence of parts-species. The **multitude of the parts** of a Platonic Idea is thus **infinite**.

*(2) Anthyphairetic Periodicity and the finite number of Logoi*

A **Logos** of a Platonic Idea is **the ratio of two successive species-parts**, namely either the ratio of the opposite species at some level, say k, $B_k/A_k$, or a ratio of the form $A_k/B_{k+1}$. Since in a Platonic Idea the Logos Criterion and the resulting periodicity hold true, there are natural numbers m>n, such that $B_n/A_n=B_m/A_m$. Thus, **the multitude of Logoi** in a Platonic Idea is **finite** and consists of all the Logoi in a period, namely:

$B_n/A_n, A_n/B_{n+1}, B_{n+1}/A_{n+1}, \ldots, B_{m-1}/A_{m-1}, A_{m-1}/B_m$ (and $B_m/A_m= B_n/A_n$).

---

[36] S. Negrepontis, *The anthyphairetic interpretation of the second hypothesis in the Parmenides*, unpublished manuscript, 2004.

The finiteness of the number of Logoi is clearly confirmed in *Sophistes* 257a4-6.

*(3) Equalisation of the Logoi by Periodicity*

The anthyphairesis of each pair of ssuccessive parts forming a Logos, say $B_k/A_k$ or $A_k/B_{k+1}$ (in the notation of (2)), is, by the philosophic analogue of the definition of the proportion of magnitudes, given in the *Topics,* periodicity, and the resulting Logos Criterion, the same, except that it appears in a cyclic permutation of the original dyad. In this way, we regard all the Logoi of a Platonic Being same and equalised to each other.

*(4) Equalization of the parts by their participation in the Logoi*

The fundamental property of a Platonic Being is that the two opposite parts—being a and non-being b—of the indefinite dyad possess the 'power' (247d-248d), which is clarified as the 'power to communicate' with each other (252d3, 253a8, 253e1, 254b8, 254c5) and still further clarified to be a 'power for equalisation with each other' (257b).
The same basic property of a Being, is described in the *Politicus* 305e-311c with regard to the Being or Form of Virtue, where the dyad a bravery and b self-restraint is called the Logos that defines Virtue.
The way the two opposing parts of the indefinite dyad of a Being (and in fact any two successive remainders in the sequence
$a, b, a_1, b_1, \ldots, a_k, b_k, a_{k+1}, \ldots$
of the anthyphairetic division) become **equalised** by Logos is well-explained in the *Sophistes* 257b1-260b3, especially in 257b, 257e9-258c5:
it is by considering
not the opposing parts themselves (which are 'not-beings', 256d-e, 257a),
but the Logoi-ratios into which they **participate**;
thus, we consider
instead of the non-being a**,** the Logos a/b, and
instead of the non-being b, the Logos a/b (or $b/a_1$) (257e2-258b5).
and we say that the part a participates into the Logos a/b, and the part b participates in the Logos a/b (oe $b/a_1$).
Thus each part participates in a Logos, and since the Logoi have been equalised by periodicity, the parts are equalised by their participation into the Logoi.

*(5) Collection of the Many into One*

Now we can understand the meaning of the Colletion of the (infinitely) many into One. Anthyphairetic Division produces, generates an infinite multitude of parts (as stated in (1). By (4) each part participates in a Logos. By periodicity (2) there is a finite number of Logoi in cyclical permutation. By (3) all these Logoi have the same anthyphairesis except for cyclical permutatyion, and in this sense they are equalised. Thus by participation of parts in Logoi all parts are equalised. This is precisely the meaning of Platonic Collection of the Many into One.
Thus, in a Platonic Idea, there is
an infinite number of parts or 'not-beings' (as confirmed in 256e5-7),

but all parts are equalised, by participation and periodicity, in the sense of self-similarity.[37]

## 14. 5. Platonic Numbers

The self-similarity of a Platonic Idea opens up the way for the definition of Platonic numbers: the anthyphairetic interpretation of Platonic numbers is a natural by-product of the anthyphairetic interpretation of the method of Division and Collection. We shall briefly outline the nature of the Numbers in accordance with the anthyphairetic interpretation.

These are certainly the same numbers described in *Philebus* 56c10-e6 (where they are differentiated from the numbers used by the many) and the *Politeia* 522c1-526c7 (where they are differentiated from the numbers as used by, say, Agamemnon). The differentiation from common usage is twofold: a Platonic number consists of units that are (i) entirely equal to each other, and (ii) divisible into an infinite multitude of parts.

The nature of these numbers is well described by Aristotle in *Metaphysics* 987b25-988a1, where it is stated that:

(a) 'the number two consists of the two elements of an indefinite dyad, equalized ('isasthenton') by the principle of the One' (1081a23-25, 1083b30-32, 1091a24-25),

(b) contrary to the mathematical units, the units of the Platonic ('eidetic') numbers are equal for a fixed number, but different for different eidetic numbers (1080a23-30), and

(c) the species are numbers (987b20-22).

The difficulties modern Platonists have encountered in reconstructing the Platonic numbers according to Aristotle's requirements has led to a downgrading of Aristotle's account of Platonic numbers (misunderstanding on his part, unwritten dogmas, and so on). But Plato describes Platonic numbers in the second hypothesis of the *Parmenides*, and the anthyphairetic interpretation of this hypothesis, discovered by S. Negrepontis[38], confirms Aristotle's description of Platonic numbers. Thus, Plato's definition of **the number two** in the *Parmenides* 144b-e as the pair of One and Being fully agrees with (a), since the self-similar Oneness of the Platonic Being precisely assures that these two successive units-species are equalized ('exisousthon', 144e2) by periodicity and Logos (or, as Aristotle calls it, the principle of the One). The units for forming a Platonic number are the species, numbered according to their successive generation in the anthyphairetic Division of the Being, precisely as in the *Parmenides* 144b-e, the *Philebus* 16c5-17a5 and the *Sophistes* 258c3, where each not-being part-species is claimed to be 'one species among the many Beings possessing number ('enarithmon'). Platonic numbers thus indeed are species, as Aristotle claims (c). The numbers in the Idea Sophist are the numbered species as in Table 14.5.

---

[37] Palmer [Pa, 1999] develops an interpretation of the Second Hypothesis of the *Parmenides*, according to which the infinite divisibility of the One Being is not an obstacle, but rather an asset, to its nature as a One, thanks to 'the demonstration of the one being's resistance to the type of division envisaged here via a demonstration of how its predicate parts "one" and "being" persist through every possible division' (p. 225). He thus approaches the self-similar Oneness of a Platonic Being, but taking into account only the infinite divisibility (which cannot by itself produce self-similarity) and not periodicity, the only cause of Platonic self-similarity.

[38] S. Negrepontis, *The anthyphairetic interpretation of the second hypothesis in the Parmenides*, unpublished manuscript, 2004.

Table 14.5. The Platonic Numbers in the Platonic Idea Sophist

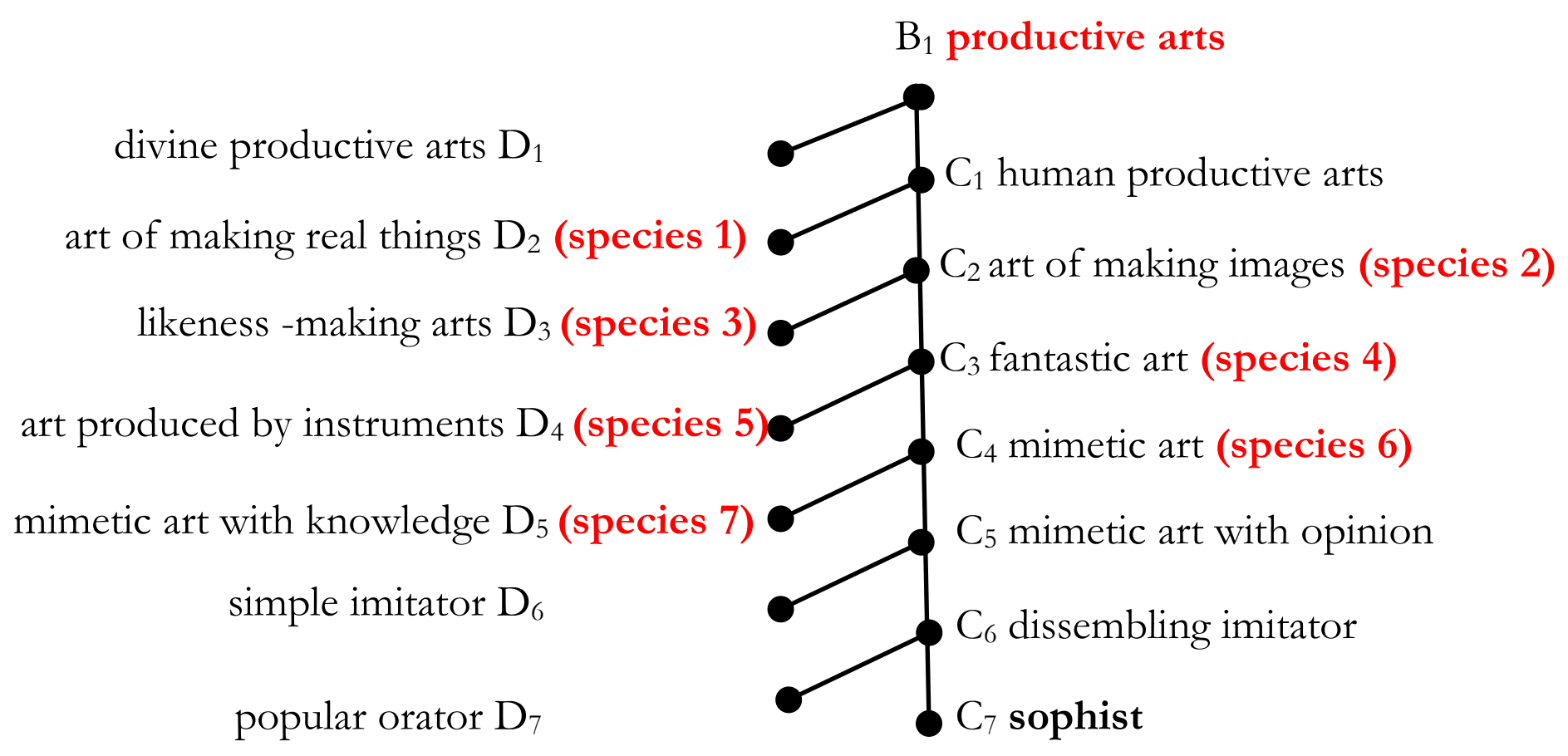

According to this interpretation, Platonic numbers are not absolute but only relative to a Platonic Being: the initial unit (hence (b)). In fact, however, according to the *Parmenides* 148d5-149d7 and the *Sophistes* 257a4-6 (‘hosaper…tosauta’, like the ‘tosauta hosaper’ in the *Parmenides* 144d5, both following the ‘tosauta…hosa’ in Zeno’s *Fragment B3*), the number of units in any Platonic Being is **finite** and equal to (the number of different Logoi)+1, this number being determined by the length of the period. In the case of the Sophist (Table 14.5) the different Logoi are D2/C2, C2/D3, D3/C3, C3/D4, D4/C4, C4/D5, while D5/C5 is equal to D2/C2, thus the number generated in the Platonic Idea of the Sophist is the number Seven.

## 14.6. Dichotomic division in Division and Collection

The Division steps in the *Sophistes* and the *Politicus* are almost[39] invariably described as dichotomic divisions[40] clearly meaning that the division is in two equal parts, not simply in two, as shown by the insistence for division in half (*Sophistes* 221b3), ‘mesotomein’ (*Politicus* 265a4), ‘dia meson’ division (*Politicus* 262b6), while on the other hand the two species into which the given genus is being defined in no way appear equal (e.g. the division of all the roductive arts into divine productive arts and human roductive arts). The answer is that these parts are equalised by the self-similarity explained above (14.5). Thus the statement in the *Politicus* that divisions should be perforned only by ‘mesotomein’ means that that the division should be performed in such a way as to ensure eventually periodicity, by means of the Logos Criterion, and thus ensure self-similarity and equalisation.

[39] Every division step in the Division and Collection is a binary division. Some scholars think that this is no longer true in the *Philebus* 16d3-5 ‘εἰ δὲ μή, τρεῖς ἤ τινα ἄλλον ἀριθμόν καὶ τῶν ἓν ἐκείνων ἕκαστον πάλιν ὡσαύτως’, but the passage refers to the number of units under the species, obtained by the *Parmenides* rule ‘number of units=(the number of different Logoi)+1’ (mentioned in 14.5), a number that, as stated in the *Philebus* 16d3-5 assage, is the same for all species. The division is described, almost invariably in every step not only as a binary division but as a dichotomic one.
[40] e.g. *Sophistes* 221e1-3, 225a2-7, 265a10-b3, 265d5-266a3, 267a1-4; *Politicus* 276d8-10, 276e6-8, 282c6-9, 287b10-c6, 302c4-11

Address of author:
*S. Negrepontis*
*Professor Emeritus,*
*Department of Mathematics,*
*Athens University, Athens 157 84, Greece*
snegrep@math.uoa.gr

May 16, 2014